\documentclass[notitlepage]{scrartcl}
\usepackage[shortlabels]{enumitem}
\usepackage{amsmath}
\usepackage{amssymb}
\usepackage{amsthm}
\usepackage{abstract}
\usepackage{xcolor}
\usepackage{comment}
\usepackage{mathtools}
\usepackage{graphicx}
\usepackage{caption}
\usepackage{subcaption}
\usepackage{float}
\usepackage{hyperref}
\usepackage{cleveref}
\usepackage{tikz}
\hypersetup{
    colorlinks=true,
    linkcolor=blue,
    filecolor=magenta,      
    urlcolor=cyan
    }
\makeatletter
\newcommand*{\barfix}[2][.175ex]{%
  \mathpalette{\@barfix{#1}}{#2}%
}
\newcommand*{\@barfix}[3]{%
  \vbox{%
    \kern#1\relax
    \hbox{$#2#3\m@th$}%
  }%
}
\makeatother

\newtheorem{theorem}{Theorem}
\newtheorem{thm}{Theorem}[section]
\newtheorem{corollary}{Corollary}
\newtheorem{lemma}[thm]{Lemma}

\newtheorem{claim}[thm]{Claim}
\newtheorem{remark}[thm]{Remark}
\newtheorem{conjecture}[thm]{Conjecture}
\newtheorem{question}[thm]{Question}

\usepackage[
    left=2cm,
    right=2cm,
    top=2cm,
    bottom=2cm,
    footskip=1cm
]{geometry}

\providecommand\given{\nonscript\:\ifthenelse{\equal{\delimsize}{}}{\big\vert}{\delimsize\vert}\nonscript\:\mathopen{}}
\let\Pr\undefined
\DeclarePairedDelimiterXPP\Pr[1]{\mathbb{P}}(){}{#1}
\DeclarePairedDelimiterXPP\Ex[1]{\mathbb{E}}{[}{]}{}{#1}

\title{\vspace{-1.5cm} Combinatorial Properties of The Giant Under Supercritical Percolation}
\author{Sahar Diskin \footnote{Department of Mathematics, ETH Z\"urich. Email: \href{mailto:Sahardiskinmail@gmail.com}{Sahardiskinmail@gmail.com}.} \and
Mihyun Kang \footnote{Institute of Discrete Mathematics, Graz University of Technology, Austria. Email: \href{mailto:kang@math.tugraz.at}{kang@math.tugraz.at}.} \and
Michael Krivelevich \footnote{School of Mathematical Sciences, Tel Aviv University, Tel Aviv 6997801, Israel. \newline Email:
\href{mailto:krivelev@tauex.tau.ac.il}{krivelev@tauex.tau.ac.il}.} \and 
Itay Markbreit \footnote{School of Mathematical Sciences, Tel Aviv University, Tel Aviv 6997801, Israel. \newline Email:
\href{mailto:markbreit@mail.tau.ac.il}{markbreit@mail.tau.ac.il}.}}

\begin{document}
\maketitle
\begin{abstract}
    We study combinatorial properties, such as diameter and long paths and cycles, of the giant component of percolated graphs. Such combinatorial properties are well studied for classical graphs such as pseudo-random graphs and the hypercube. We aim to obtain similar results for extensive families of $d$-regular $n$-vertex graphs. Specifically, assuming expansion properties of the host graph, we show a typical vertex expansion by a factor of $1/d\log d$ in the giant component of the percolated subgraph. As a consequence, we obtain bounds on various combinatorial properties of the giant component. For example, a bound of $O\left(d\log^2 d+\log n\right)$ on the typical diameter of the giant component of the percolated subgraph. In addition, we construct graphs with very good vertex and edge expansion, but having large diameter, showing that the diameter of the percolated subgraph is not necessarily typically logarithmic in the number of the vertices of the graph. We also present good expanders such that their percolated subgraph typically contains a bad expanding set.
\end{abstract}

\section{Introduction}\label{intro}

\subsection{Background and Motivation}

In 1957, Broadbent and Hammersley \cite{zbMATH03148802} initiated the study of percolation theory in order to model the flow of fluid through a medium with randomly blocked channels. Since then, the theory of percolation has been studied extensively (see \cite{zbMATH05076923}, \cite{zbMATH01301486} and \cite{zbMATH03826957} for systematic coverage).

Given a host graph $G=(V,E)$ and probability $p\in [0,1]$, in \textit{$p$-bond percolation} we form a random subgraph $G_p\subseteq G$ by retaining each edge of $G$ independently with probability $p$. We call $G_p$ the percolated graph.

The most well-studied example of the percolated graph is the binomial random graph $G_{n,p}$ when the host graph is the complete graph $K_n$. A classical result of Erd\H os and R\'enyi from 1960 \cite{zbMATH03168330} states that $G_{n,p}$ undergoes a fundamental phase transition with respect to its connected component structure when the expected average degree is around $1$ (that is, $n\cdot p \approx 1$): for any constant $\epsilon>0$, when $p=\frac{1-\epsilon}{n}$, all the components are typically of order $O\left(\epsilon^{-2}\log n\right)$, but when $p=\frac{1+\epsilon}{n}$, there typically exists a unique giant component $L_1$ in $G_{n,p}$ of order $(1 + o(1))\gamma n$, and all other components of $G_{n,p}$ are of order $O\left(\epsilon^{-2}\log n\right)$, where
\begin{align}\label{gamma_def}
    \gamma\coloneqq \gamma(\epsilon)=2\epsilon-O\left(\epsilon^2\right)
\end{align}
is the survival probability of the $\mathrm{Po}(1+\epsilon)$ branching process. For more background on branching processes, see, for example, \cite{zbMATH03410334} and \cite[Section 4.3.4]{zbMATH06976511}.

Since then, the phase transition discovered by Erd\H os and R\'enyi has been established for several other families of graphs, including pseudo-random $(n,d,\lambda)$-graphs \cite{zbMATH02037533} and the $d$-dimensional binary hypercube $Q^d$ \cite{zbMATH03769676,zbMATH00035138,arXiv:2311.07210}. More precisely, if $p$-bond percolation is performed with $p=\frac{1+\epsilon}{d}$, where $\epsilon>0$ is a sufficiently small constant and $d$ is a sufficiently large integer, then, \textbf{whp}\footnote{With high probability, that is, with probability tending to one as $n$ tends to infinity.}, the percolated graph contains a unique giant component, where its size, and the size of the remaining components match the sizes of the corresponding supercritical random graph.

More recently, Diskin and Krivelevich \cite{zbMATH08242616} identified a general set of conditions on the host graph $G$ guaranteeing the same conclusion on the percolated graph $G_p$. In addition, Diskin, Erde, Kang, and Krivelevich \cite{zbMATH08225676} discovered another general set of properties of the host graph $G$, under which \textbf{whp} $G_p$ contains a unique giant component of the correct order and all remaining components are substantially smaller than the giant component. These results generalized the previous results where the host graph is $K_n$, an $(n,d,\lambda)$-graph, or $Q^d$. We give here the exact statements from \cite{zbMATH08225676,zbMATH08242616}, using some standard notations defined at the beginning of Section \ref{pre}:

\begin{thm}[Theorem 1 of \cite{zbMATH08225676}]\label{isop}
    Let $d = \omega_n(1)$, and let $G$ be a $d$-regular $n$-vertex graph. Let $\epsilon > 0$ be a small constant and let $p=\frac{1+\epsilon}{d}$. Then, there exists a constant $C_0 \coloneqq C_0(\epsilon) > 0$, such that for every $C \ge C_0$ the following holds. If for every $U\subseteq V(G)$ with $|U|\le \frac{n}{2}$, we have $e_G(U,U^c)\ge C|U|$, then with probability at least $1-o_d(1)-\exp\left(-\frac{\sqrt{C}n}{d}\right)$, the subgraph $G_p$ contains a giant component $L_1$, satisfying 
    \begin{align*}
        \Big||V(L_1)|-\gamma\cdot n\Big| \le \frac{10n}{C^{1/20}},
    \end{align*}
    where $\gamma\coloneqq \gamma(\epsilon)$ is as defined in \eqref{gamma_def}. Furthermore, all other components in $G_p$ are of size at most $(1+o_d(1))\frac{10n}{C^{1/20}}$.
\end{thm}

\begin{thm}[Theorem 1 of \cite{zbMATH08242616}]\label{ERCP1}
    Let $\alpha,c_1,c_2,c_3>0$ be constants and let $\epsilon>0$ be a sufficiently small constant. Let $d,n$ be sufficiently large integers such that $d\ge \log^\alpha n$. Let $G$ be a $d$-regular $n$-vertex graph, and let $p=\frac{1+\epsilon}{d}$. Suppose there exists a sufficiently large $C\coloneqq C(\epsilon,c_1,c_2,c_3,\alpha)>0$ that $G$ satisfies the following vertex and edge expansion properties.
    \begin{enumerate}[(P\arabic*{})]
        \item For every $U\subseteq V(G)$ with $|U|\le \frac{n}{2}$, $e_G(U,U^c)\ge c_1|U|$; \label{prop0P}
        \item For every $U\subseteq V(G)$ with $|U|\le c_2\log n$, $|N_G(U)|\ge c_3d|U|$; \label{prop1P}
        \item For every $U\subseteq V(G)$ with $|U|\le Cd\log n$, $e_G(U,U^c)\ge (1-\epsilon^3)d|U|$. \label{prop2P}
    \end{enumerate}
    Then, \textbf{whp} there exists a unique giant component $L_1$ in $G_p$ whose size is $(1 + o(1))\gamma n$, where $\gamma\coloneqq \gamma(\epsilon)$ is as defined in \eqref{gamma_def}, and all other components in $G_p$ are of size $O\left(\frac{\log n}{\epsilon^2}\right)$.
\end{thm}

\begin{remark}\label{re: sprinkling_ercp1}
    Although property \ref{prop2P} in Theorem~\ref{ERCP1} is stated with $1-\epsilon^3$, the proof of Theorem~\ref{ERCP1} in \cite{zbMATH08242616} actually relies only on the order of magnitude of $\epsilon^3$. Consequently, the exact same argument holds under the slightly weaker assumption with $1-2\epsilon^3$, as the precise constant before $\epsilon^3$ is flexible within their proof.
    
    Later in this paper, we consider graphs satisfying properties \ref{prop0P}-\ref{prop2P} and apply Theorem~\ref{ERCP1} with $p=\frac{1+\epsilon-\delta}{d}$, where $\delta\coloneqq \delta(\epsilon)>0$ is a small constant in terms of $\epsilon$. The aforementioned mild relaxation of property \ref{prop2P} to $1-2\epsilon^3$ allows us to apply Theorem~\ref{ERCP1} in our setting.
\end{remark}

Theorem~\ref{isop} provides, in a sense, a relatively weak sufficient condition on the host graph $G$ that guarantees the typical emergence of a unique giant component in the random subgraph $G_p$. Our first aim in this paper is to establish another relatively weak set of conditions on a $d$-regular $n$-vertex host graph $G$ (where $d=\omega_n(1)$), which are also tight and yield essentially the same conclusion as in Theorem~\ref{isop} on the percolated graph $G_p$.

Theorem~\ref{ERCP1} in turn gives a full characterization of the typical sizes of the components of the percolated graph under supercritical percolation, establishing \textbf{whp} the unique giant component of the `correct' order of $\gamma n$, as well as a large gap between the size of the giant component and the size of all other components, which is  $O\left(\epsilon^{-2}\log n\right)$, achieving essentially the same phenomenon Erd\H os and R\'enyi discovered for $G_{n,p}$. Consequently, naturally this phenomenon is sometimes referred to as ERCP --- Erd\H os-R\'enyi component phenomenon. 

The main goal of our paper is to uncover the typical rich combinatorial properties of the giant component in this setting. We begin by establishing relatively strong vertex expansion, which in turn enables us to derive a range of further interesting combinatorial aspects, including the existence of long cycles, large clique minors, and a small diameter. We then investigate the tightness of these results and establish matching statements showing that our bounds are, in several respects, close to optimal. These results also show that the settings considered in Theorem~\ref{L1} and Theorem~\ref{isop}, which do not establish the full ERCP, are not sufficient to guarantee \textbf{whp} the rich combinatorial properties of the giant component that hold under the setting of Theorem~\ref{ERCP1}.

The inspiration for this study comes from several known results asserting rich combinatorial properties and aspects of the giant component of the percolated graph for well-known host graphs such as $K_n$, an $(n,d,\lambda)$-graph, or $Q^d$, where $n$ and $d\coloneqq d(n)$ are sufficiently large integers. Let us elaborate on the known results for each host graph individually:

\paragraph{The complete graph $K_n$.} The $p$-bond percolated $K_n$ is the binomial random graph $G_{n,p}$. Chung and Lu \cite{zbMATH01623567}, Fernholz and Ramachandran \cite{zbMATH05224135}, and Riordan and Wormald \cite{zbMATH06144670} showed, independently, that \textbf{whp} the diameter of the giant component of $G_{n,p}$ is $\Theta_\epsilon(\log n)$. Fountoulakis and Reed \cite{zbMATH05320778}, and independently Benjamini, Kozma, and Wormald \cite{zbMATH06370216} obtained a good control over the typical edge-expansion of connected subsets of the giant component in $G_{n,p}$ and utilized it to show that the asymptotic typical order of magnitude of the mixing time of a simple random walk on the giant component is asymptotically $\Theta_\epsilon\left(\log^2 n\right)$.

\paragraph{Pseudo-Random $(n,d,\lambda)$-graphs.}
Let us first recall the notion of $(n,d,\lambda)$-graphs. Assume that the eigenvalues of (the adjacency matrix of) a graph $G$ are ordered in the non-increasing order $\lambda_1 \geq \dots \geq \lambda_n$. The largest eigenvalue of any $d$-regular graph is easily seen to be $d$, sometimes referred to as the trivial eigenvalue of $G$. A graph $G$ is an $(n, d, \lambda)$-\textit{graph} if $G$ is a $d$-regular $n$-vertex graph, and all the eigenvalues of $G$, but the first one, are at most $\lambda$ in their absolute values. We refer the reader to the survey \cite{zbMATH05037081} for an extensive discussion on $(n,d,\lambda)$-graphs and other notions of pseudo-random graphs. Diskin and Krivelevich \cite{zbMATH07931146} established an asymptotic typical vertex expansion for sets inside the giant component of the $p$-bond percolated $(n,d,\lambda)$-graph. They showed that typically sets of linear size, as well as connected sets of size between logarithmic in $n$ and linear in $n$ have vertex expansion factor of $\tilde{\Omega}\left(\epsilon^2\right)$ inside the giant component. From these expansion properties, they derive that \textbf{whp} the diameter of the giant component is $O_\epsilon(\log n)$, and that the typical mixing time of a lazy simple random walk on the giant component is asymptotically $O_\epsilon(\log^2 n)$. 

\paragraph{The $d$-dimensional binary hypercube $Q^d$.}
The $d$-dimensional hypercube $Q^d$ is the graph with the vertex set $V(Q^d)=\{0,1\}^d$, where two vertices are adjacent if they differ in exactly one coordinate. Erde, Kang, and Krivelevich \cite{zbMATH07628800} showed that \textbf{whp} the vertex expansion of the giant component of $Q^d_p$ is inverse polynomial in $d$. As a consequence, they obtained \textbf{whp} polynomial in $d$ bounds on the diameter of the giant component and the mixing time of a lazy simple random walk on the giant component, as well as lower bounds on the circumference and Hadwiger number of the giant component. See also \cite{zbMATH07915715} which strengthen these results and extend them to more general product graphs. Finally, Anastos, Diskin, Lichev, and Zhukovskii \cite{arXiv:2510.13348} improved these results and resolved the long-standing open problem regarding the typical diameter and the typical mixing time of a lazy simple walk on the giant component of $Q^d_p$. They showed that \textbf{whp} the diameter of the giant component is $\Theta(d)$ and that the mixing time of the lazy simple random walk on the giant component is $\Theta\left(d^2\right)$.

\subsection{Main Results} \label{main results}

Throughout the paper, we let $\epsilon>0$ be a sufficiently small constant. Let $G$ be a $d$-regular $n$-vertex graph, where $n,d\in \mathbb{N}$ are sufficiently large and satisfy $d = \omega_n(1)$ and $d=o(n)$. Set $p\in[0,1]$. Denote by $L_1\coloneqq L_1(G_p)$ the largest component of $G_p$. In this section we further set $p=\frac{1+\epsilon}{d}$.

In Theorems \ref{isop} and \ref{ERCP1}, as well as in our results, we assume that the host graph $G$ satisfies vertex and/or edge expansion properties. To simplify the presentation, we introduce notation for the assumptions considered throughout the paper:

\begin{enumerate}
    \item [(A1)] Let $b>0$. We say that $G$ satisfies \textbf{iso}$\mathbf{(b)}$ if for every $U\subseteq V(G)$ with $|U|\le \frac{n}{2}$, $e_G(U,U^c)\ge b|U|$.
    \item [(A2)] Let $A,\alpha>0$. We say that $G$ satisfies \textbf{VE}$\mathbf{(A,\alpha)}$ if for every $U\subseteq V(G)$ with $|U|\le A$, $|N_G(U)|\ge (1-\alpha)d|U|$.
    \item [(A3)] Let $A,\alpha>0$. We say that $G$ satisfies \textbf{EE}$\mathbf{(A,\alpha)}$ if for every $U\subseteq V(G)$ with $|U|\le A$, $e_G(U,U^c)\ge (1-\alpha)d|U|$.
\end{enumerate}

Using the notation above, let us note that properties \ref{prop0P}-\ref{prop2P} in Theorem~\ref{ERCP1} can be reformulated as $\mathrm{iso}(c_1)$, $\mathrm{VE}(c_2\log n,1-c_3)$ and $\mathrm{EE}(Cd\log n,\epsilon^3)$.

Our first result provides a third set of conditions on the host graph which guarantees the typical appearance of a unique giant component of the percolated graph.

\begin{theorem}\label{L1}
    Let $C,b>0$ such that $Cb$ is sufficiently large in terms of $\epsilon$. Set $0<\alpha\le \epsilon/4$. Let $G$ be a $d$-regular $n$-vertex graph satisfying $\mathrm{iso}(b)$ and $\mathrm{EE}(Cd,\alpha)$. Then, \textbf{whp} the giant component $L_1$ of $G_p$ satisfies
    \begin{align*}
        \Big||V(L_1)|-\gamma\cdot n\Big| \le \frac{10}{(Cb)^{1/20}} n,
    \end{align*}
    where $\gamma\coloneqq \gamma(\epsilon)$is as defined in \eqref{gamma_def}. Furthermore, all other components in $G_p$ are of size at most $(1+o_d(1))\frac{10}{(Cb)^{1/20}} n$.
\end{theorem}

A few comments are in place.  The upper bound $\epsilon/4$ on the parameter $\alpha$ in Theorem~\ref{L1} is not tight and can be increased, as long as $\alpha$ remains sufficiently far from $\epsilon$. In Theorem~\ref{L1} we allow $C$ or $b$ to depend on $d,n$. In particular, for example when $C=C(d,n)\to \infty$, the
result implies that \textbf{whp} $|L_1| = (1 + o_d(1)) \gamma n$ and all other components have sublinear order. We further note that the framework of Theorem~\ref{L1} extends that of Theorem~\ref{isop}. Indeed, Theorem~\ref{L1} applies for any positive constant $b$, whereas Theorem~\ref{isop} requires $b$ to be sufficiently large in terms of $\epsilon$. Consequently, Theorem~\ref{L1} applies immediately to a wider class of graphs, including the $d$-dimensional binary hypercube $Q^d$ (which does not satisfy the assumptions of Theorem~\ref{isop}, recalling that $\mathrm{iso}(Q^d)= 1$). Nevertheless, there is a tradeoff between the two sets of assumptions on the host graph $G$. In Theorem~\ref{L1}, $b$ can be any constant, allowing $b$ to be small, but this comes at the cost of requiring the constant $C$ to be larger. Finally, we note that Theorem~\ref{L1} is tight up to the choice of the parameter $C$ --- the full statement can be found in Section \ref{tightness}.

Next we examine the strength of the settings presented in Theorems \ref{isop}, \ref{ERCP1}, and \ref{L1} in establishing typical combinatorial properties of the giant component of the percolated graph. The following results in this section show that the setting of Theorem~\ref{ERCP1} suffices to guarantee lower bounds on the typical vertex expansion and upper bounds on the typical diameter of the giant component of the percolated graph. Conversely, the results in Section \ref{tightness} provide constructions, relevant to all three settings of Theorems \ref{isop}, \ref{ERCP1}, and \ref{L1}, demonstrating the limits of obtaining strong vertex expansion and a small diameter of the giant component, as well as showing the relative tightness of our results under the setting of Theorem~\ref{ERCP1}.

The following result provides the typical vertex expansion of the giant component of the percolated graph, when its host graph satisfies the assumptions of Theorem~\ref{ERCP1}.

\begin{theorem}\label{general}
    Under the notations and assumptions of Theorem~\ref{ERCP1}, there exists a constant $\beta\coloneqq \beta(\epsilon,c_1,c_2,c_3,\alpha)>0$ such that \textbf{whp} for every $S\subseteq V(L_1)$ with $\epsilon^2n\le |S|\le \frac{3\epsilon n}{2}$,
    \begin{align*}
        |N_{G_p}(S)|\ge \frac{\beta}{d\log d}|S|.
    \end{align*}
\end{theorem}

The following corollary is a (by now) standard consequence from Theorem~\ref{general}:

\begin{corollary}\label{full expansion in subgraph}
    Under the notations and assumptions of Theorem~\ref{general}, there exists a constant $\beta\coloneqq \beta(\epsilon,c_1,c_2,c_3,\alpha)>0$ such that \textbf{whp} there exists a subgraph of the giant component $H\subseteq L_1$ with $|V(H)| \ge \frac{3\epsilon n}{2}$ such that for every $S \subseteq V(H)$ with $|S|\le \frac{|V(H)|}{2}$,
    \begin{align*}
        |N_{H}(S)|\ge \frac{\beta}{d\log d}|S|.
    \end{align*}
\end{corollary}

Moreover, the vertex expansion property we obtained in Theorem~\ref{general} has several interesting standard consequences for combinatorial properties of the giant component $L_1$.

\begin{corollary}\label{properties_general}
    Under the notations and assumptions of Theorem~\ref{general}, \textbf{whp}:
    \begin{enumerate}[(\arabic*)]
        \item The giant component $L_1$ of $G_p$ contains a cycle of length $\Omega\left(\frac{n}{d\log d}\right)$;
        \item \label{cor_2} The Hadwiger number of $L_1$ is 
        \begin{align*}
            \Omega\left(\frac{\sqrt{n}}{d\log d}\right).
        \end{align*}
    \end{enumerate}
\end{corollary}

We note that \textbf{whp} $|E(G_p)|= O(n)$, and the Hadwiger number $h(H)$ of a graph $H$ satisfies $|E(H)|\ge \binom{h(H)}{2}$. Hence, both results in Corollary~\ref{properties_general} are optimal up to
polynomial terms in $d$.

Our next theorem establishes typical vertex expansion properties of connected subgraphs of the giant component of the percolated graph under the assumptions in Theorem~\ref{ERCP1}.

\begin{theorem}\label{setup_2_expansion}
    Under the notations and assumptions of Theorem~\ref{ERCP1}, and for every constant $C_1>0$, there exist constants $C_2\coloneqq C_2(C_1,c_1,c_2,c_3,\alpha,\epsilon)>0$ and $\beta\coloneqq \beta(C_1,c_1,c_2,c_3,\alpha,\epsilon)>0$ such that \textbf{whp} for every subset of the giant component $S\subseteq V(L_1)$ with $nd^{-C_1}\le |S|\le \frac{3\epsilon n}{2}$ and $|S|\ge C_2 (d\log d+\log n)$ and such that $G_p[S]$ is connected,
    \begin{align*}
        |N_{G_p}(S)|\ge \frac{\beta}{d\log d}|S|.
    \end{align*}
\end{theorem}

A few comments are in place. Note that in Theorem~\ref{setup_2_expansion}, $C_1$ can be chosen to be any constant, with $\beta$ depending on this choice. In fact, the proof of Theorem~\ref{setup_2_expansion} allows us to take $C_1$ as large as $\log\log d$, thereby allowing a wider range of sets $S$ in the statement, at the cost of weakening the vertex expansion factor from $\beta/d\log d$ to $\beta'/d\log d\log\log d$ with $\beta'\coloneqq \beta'(c_1,c_2,c_3,\alpha,\epsilon)$ independent of $d$. It turns out that Theorem~\ref{setup_2_expansion} is tight under similar assumption on the host graph $G$ as can be seen below in Theorem~\ref{no_expansion}. The following result uses Theorem~\ref{setup_2_expansion} to guarantee an upper bound on the typical diameter of the giant component of the percolated graph.

\begin{corollary}\label{positive_diam}
    Under the notations and assumptions of Theorem~\ref{setup_2_expansion}, and if we further assume that $G$ satisfies $\mathrm{EE}\left(nd^{-C_1},\epsilon^3\right)$, then, \textbf{whp}
    \[
        \mathrm{diam}(L_1)=O\left(d\log^2 d+\log n\right).
    \]
    
\end{corollary}

It should be noted that the estimate for the typical diameter of $L_1$ in Corollary~\ref{positive_diam} is tight up to a multiplicative constant that depends on $\epsilon$ and on $\beta$.

Let us present another result concerning upper bound on the typical diameter of the giant component of the percolated graph under the setting of Theorem~\ref{ERCP1}.

\begin{theorem}\label{pure_diam}
    Under the notations and assumptions of Theorem~\ref{ERCP1}, \textbf{whp} 
    \begin{align*}
        \mathrm{diam}(L_1)=O\left(d\log n\left(1+\log_d\log n\right)\right).
    \end{align*}
\end{theorem}

Notice that Theorem~\ref{pure_diam} shows in particular that, whenever $d=\Omega(\log^c n)$ for some (small) constant $c>0$, the diameter of the giant component $L_1$ satisfies \textbf{whp} $\mathrm{diam}(L_1)=O(d\log n)$.

\paragraph{Proof strategy.} Let us briefly explain the main idea behind the proofs of the four theorems stated above. All proofs of these results use a sprinkling argument, where the first exposure is still supercritical. Roughly speaking, this allows us to obtain almost everything we need after the first exposure, while the second exposure serves as a refinement that yields the full results. Throughout the proofs, we combine several techniques, allowing us to obtain robust results for a broad class of graphs.

At a high level, in Theorem~\ref{L1}, the first exposure ensures that \textbf{whp} the vast majority of vertices in the host graph are close (in graph distance) to many vertices that lie in `large' components of the percolated graph. The second exposure merges these `large' components into a unique giant one. 
In Theorems \ref{general} and \ref{setup_2_expansion} we show that sets in the giant component obtained after the first exposure, which we call the \textit{early giant}, already typically have good vertex expansion (thanks to edges from both exposures), then use it to obtain typical good vertex expansion for a general set in the giant component obtained after the second exposure. 
In Theorem~\ref{pure_diam}, we consider the early giant obtained after the first exposure and start with a vertex in it. We use expansion properties of the host graph to accumulate \textbf{whp} a large set of vertices that are close to the starting vertex via edges from both exposures, and iterate this procedure until at least half of the early giant is covered, thus bounding from above the typical diameter of the early giant (using edges from both exposures). Finally, we show that typically every vertex in the giant component after the second exposure is close to a vertex in the early giant from the first exposure, implying the desired upper bound on the diameter.

\subsection{Tightness and Counterexamples}\label{tightness}

In this section we discuss the tightness of the theorems presented in the previous section. Specifically, we provide constructions across the main three studied areas from Section \ref{main results}: the vertex expansion of the giant component of the percolated graph, its diameter, and the parameter regimes under which a giant component is completely prevented from forming \textbf{whp}.

We start by showing the tightness of Theorem~\ref{L1}. More concretely, we show that there exists a $d$-regular $n$-vertex graph $G$ that satisfies the assumptions of Theorem~\ref{L1}, but with a sufficiently smaller $C>0$ compared to the choice in Theorem~\ref{L1}, and \textbf{whp} $G_p$ contains no linear-sized component in the supercritical regime, and in fact all components of $G_p$ are of size $O(d\log n)$, for some supercritical probability $p$.

\begin{theorem}\label{no_L1}
    Let $\sigma>1,\alpha>0$ be some constants. Then, there exist constants $0<b<\frac{1}{\sigma}$ and $1<C=\Theta\left(\frac{\alpha}{\sigma b}\right)$ such that the following holds. Let $d = \omega_n(1)$, $d=o\left(\sqrt{n}\right)$ be an even integer and let $p=\frac{\sigma}{d}$. There exists a $d$-regular $n$-vertex graph $G$ that satisfies $\mathrm{iso}(b)$ and $\mathrm{EE}(Cd,\alpha)$, and, moreover, \textbf{whp} all components of $G_p$ are of size at most $\frac{9}{\alpha}Cd\log n$.
\end{theorem}

The following result presents the relative tightness of the typical vertex expansion factor of the giant component of the percolated graph, under some general vertex and edge expansion assumptions on the host graph that is a $d$-regular $n$-vertex graph.
Perhaps somewhat surprisingly, in some instances of graphs under assumptions that are slightly weaker than those in Theorem~\ref{ERCP1}, we show that for $d\le n^\alpha$, for some $\alpha>0$, we can typically find a connected subset of the giant component with vertex expansion of order at most $\frac{1}{d\log n}$. In particular, the following result, combined with Theorem~\ref{setup_2_expansion} in the case where the degree of the host graph is polynomial with $n$, shows that the typical vertex expansion factor $\frac{1}{d\log d}$ of the giant component in Theorem~\ref{setup_2_expansion} is tight up to a multiplicative constant.

\begin{theorem}\label{no_expansion}
    Let $\epsilon>0$ be a sufficiently small constant. For any constants $\alpha>0$ and $0<c_1<\frac{1}{160\epsilon}$, there exists a constant $c_2\coloneqq c_2(\alpha,\epsilon)>0$ such that the following holds. Let $n,d$ be sufficiently large integers such that $dn$ is even, $d = \omega_n(1)$ and $d<\sqrt{n}$. There exists a $d$-regular $n$-vertex graph $G$ that satisfies $\mathrm{iso}(c_1)$, $\mathrm{VE}(c_2\log n,\alpha)$, and $\mathrm{EE}\left(\frac{\alpha}{2}\cdot d\log n,\alpha\right)$, and, moreover, there exists a constant $\beta\coloneqq \beta(\epsilon)\ge 1$ such that for $p=\frac{1+\epsilon}{d}$, \textbf{whp} there exists a subset of the giant component $S\subseteq V(L_1)$ with $|S|=\Theta(d\log n)$ such that $G_p[S]$ is connected and
    \begin{align*}
        |N_{G_p}(S)|\le \frac{\beta}{d\log n}|S|.
    \end{align*}
\end{theorem}

A few remarks are in order. First, let us note that the existence of a giant component in Theorem~\ref{no_expansion} is guaranteed by the conclusion of Theorem~\ref{L1}. Furthermore, observe that the assumptions of Theorem~\ref{no_expansion} and Theorem~\ref{ERCP1} are very similar, differing only in the constant factor preceding the term $d\log n$ in the edge expansion property (namely, in \ref{prop2P} and $\mathrm{EE}\left(\frac{\alpha}{2}\cdot d\log n,\alpha\right)$). In particular, the same construction with a minor modification of the proof of Theorem~\ref{no_expansion} shows that \textbf{whp} the second-largest component of the percolated graph has order at least $d\log n$.

We also note that Theorem~\ref{no_expansion} works in a more general regime, where $p=\frac{\sigma}{d}$ for any constant $\sigma>1$. In this regime, we will have that \textbf{whp} there exists a subset $S$ of a linear-sized component with $|S|=\Theta(d\log n)$ such that $G_p[S]$ is connected and
\begin{align*}
    |N_{G_p}(S)|\le \frac{\beta}{d\log n}|S|.
\end{align*} 

The following two results demonstrate the relative tightness of the diameter-related statements presented in Section \ref{main results}, under certain conditions. The first result is a general statement showing in particular that if, for a general $d$-regular $n$-vertex graph $G$, a random subgraph $G_p$ contains a linear-sized component, then its diameter is \textbf{whp} $\Omega\left(\frac{\log n}{\log \sigma}\right)$. This, in turn, shows that Corollary~\ref{positive_diam} is tight for $d=O\left(\frac{\log n}{(\log\log n)^2}\right)$.

\begin{theorem}\label{like_Chung_Lu}
    Let $\sigma>1$ be a constant. Let $n$ be a sufficiently large integer and let $3\le d\coloneqq d(n)\in \mathbb{N}$. Let $G$ be a $d$-regular $n$-vertex graph and let $p=\frac{\sigma}{d-1}$. Then, \textbf{whp} every $v\in V(G)$ satisfies
    \[
        \left|B_{G_p}\left(v,\frac{\log n}{2\log \sigma}\right)\right|=o(n).
    \]
\end{theorem}

The proof of Theorem~\ref{like_Chung_Lu} is very similar to the proof of Lemma 1 and Theorem 6 in \cite{zbMATH01623567}, so we omit it.

Our final result provides a construction of an expander (both in edges and vertices) whose largest component after $p$-bond percolation has \textbf{whp} a quite large diameter. To be precise, the typical diameter of the largest component $L_1$ will be $\Omega\left(\min \left\{\frac{d}{n},\frac{1}{d}\right\}\cdot |V(L_1)|\right)$. This, in turn, shows that Corollary~\ref{positive_diam} and Theorem~\ref{pure_diam} are relatively tight, as reflected in Corollary~\ref{large_diam_setup2}. 

Furthermore, the following result emphasizes the difference between the settings in Theorems \ref{isop}, \ref{ERCP1}, \ref{L1} on the one hand and important properties of the complete graph $K_n$ and the family of $(n,d,\lambda)$-graphs on the other hand, for which, under standard assumptions on the parameters, the diameter after percolation is typically $\Theta(\log n)$ (see, for example, \cite{zbMATH01623567,zbMATH05224135,zbMATH06144670} for $K_n$, and \cite{arXiv:2510.13348,zbMATH07931146} for the hypercube and the family of $(n,d,\lambda)$-graphs, respectively).

\begin{theorem}\label{general_large_diam}
    Let $c_1,\alpha>0$ and $\sigma>1$ be constants. Let $n,d$ be sufficiently large integers such that $d = \omega_n(1), d=o\left(n\right)$. Then, for every $p\ge \frac{\sigma}{d}$ and for every even integer $\max\left\{\frac{4c_1}{\alpha}\cdot\frac{n}{d},2d\right\}\le t< n$, there exist a constant $c_2\coloneqq c_2(\alpha)>0$ and a $d$-regular $n$-vertex graph $G$ that satisfies $\mathrm{iso}(c_1)$, $\mathrm{VE}(c_2\cdot \frac{t}{d},\alpha)$, and $\mathrm{EE}(\frac{\alpha}{2}t,\alpha)$, and, moreover, the diameter of the largest component $L_1$ in $G_p$ is at least $\frac{|V(L_1)|}{2t}-1$.
\end{theorem}

We note that the proof of Theorem~\ref{general_large_diam} shows an even stronger result. It shows that \textbf{whp} there are two vertices in $L_1$ whose distance in $G$ is at least $\frac{|V(L_1)|}{2t}-1$.

Observe that by choosing $t=C\max\left\{\frac{n}{d},d\right\}$ where $C>0$ is a constant, the assumptions of Theorem~\ref{general_large_diam} encompass those of Theorems \ref{isop} (for any fixed $C>0$) and \ref{L1} (for any sufficiently large constant $C\coloneqq C(\epsilon)>0$). Consequently, we obtain the following corollary.

\begin{corollary}\label{large_diam_setup1}
    Under the notations of Theorem~\ref{isop} or \ref{L1}, for every sufficiently large integers $n,d$ such that $d = \omega_n(1), d=o\left(n\right)$, there exists a $d$-regular $n$-vertex graph $G$, satisfying the assumptions of Theorem~\ref{isop} or \ref{L1}, and moreover, the diameter of the largest component $L_1$ in $G_p$ is \textbf{whp}
    \[
        \Omega\left(\min\left\{d,\frac{n}{d}\right\}\right).
    \]
\end{corollary}

Similarly, by choosing $t=C\max\left\{\frac{n}{d},d\log n\right\}$ where $C\coloneqq C(\epsilon)>0$ is any sufficiently large constant, the assumptions of Theorem~\ref{general_large_diam} cover those of Theorem~\ref{ERCP1}.

\begin{corollary}\label{large_diam_setup2}
    Under the notations of Theorem~\ref{ERCP1}, for every sufficiently large integers $n,d$ such that $d = \omega_n(1), d=o\left(n\right)$, there exists a $d$-regular $n$-vertex graph $G$, satisfying the assumptions of Theorem~\ref{ERCP1}, and moreover, the diameter of the largest component $L_1$ in $G_p$ is \textbf{whp}
    \[
        \Omega\left(\min\left\{d,\frac{n}{d\log n}\right\}\right).
    \]
\end{corollary}

In particular, Theorem~\ref{pure_diam} together with Corollary~\ref{large_diam_setup2} establishes tight estimate for the diameter of the giant component up to a factor of $\log n$, when the host graph satisfies the assumptions of Theorem~\ref{ERCP1}, assuming $\log^c n\le d\le \sqrt{\frac{n}{\log n}}$ for some (small) constant $c>0$.

\subsection{Main Contributions of the Paper}

In this paper, we study a variety of typical combinatorial properties of the giant component of the percolated graph for a broad family of $d$-regular $n$-vertex graphs that satisfy vertex and edge expansion properties, substantially generalizing previous results for specific well-known classes of graphs.

Specifically, we prove a nearly optimal typical vertex expansion factor of $\frac{1}{d\log d}$ for the giant component of the percolated graph. While one cannot hope for a vertex expansion factor better than $\frac{1}{d}$, we show that even a $\frac{1}{d}$ vertex expansion factor is unattainable in general. Specifically, we construct an example of a well-behaved graph that \textbf{whp} contains a connected subgraph in the giant component of the percolated graph with poor vertex expansion. 

Additionally, we establish \textbf{whp} a polynomial bound in $d$ and $\log n$ on the diameter of the giant component of the percolated graph. Furthermore, perhaps surprisingly, we construct well-behaved graphs whose diameter and the typical diameter of the giant component of the percolated graph is at least $d$ (where $d\le \sqrt{n/\log n}$), marking a fundamental departure from the typical diameter bounds observed in well-known graph families, such as the giant component of the percolated graph of $K_n$ or an $(n,d,\lambda)$-graph.

\subsection{Structure of the Paper}

The paper is structured as follows. In Section \ref{pre} we set our notation and collect several useful lemmas. In Section \ref{pL1} we prove that the third setting presented in Theorem~\ref{L1} is sufficient to ensure the appearance of a unique giant component of the percolated graph. In Section \ref{Pgeneral} we provide the proofs of Theorem~\ref{general}, Corollary~\ref{full expansion in subgraph}, and Corollary~\ref{properties_general}. Section \ref{Psetup_2_expansion} is devoted to the proof of Theorem~\ref{setup_2_expansion} and Corollary~\ref{positive_diam}. In Section \ref{Ppure_diam} we prove Theorem~\ref{pure_diam}. Sections \ref{pno_L1}, \ref{Pno_expansion} and \ref{P_larege_diam} are devoted to the tightness results regarding the absence of a giant in the percolated graph, the vertex expansion factor and the diameter of the largest component of the percolated graph (respectively). We conclude with a final discussion in Section \ref{disc}. 

\section{Preliminaries}\label{pre}

Let us first state some notational conventions used in the paper. Consider a graph $H=(V,E)$ and a vertex $v \in V$. For an integer $r$, let $B_H(v, r)\subseteq V(H)$ be the ball of radius $r$ in $H$ centered at $v$, which we consider as the set of vertices of $H$ at a distance of at most $r$ from the center $v$. We also denote by $C_H(v)$ the connected component of $H$ that contains $v$. As is fairly standard, $E_H(S, S^c)$ denotes the set of edges in $E$ with one endpoint in $S$ and the other in $S^c\coloneqq V(H) \setminus S$, and $E_H(S)$ denotes the set of edges in the induced subgraph $H[S]$. We use $e_H(S, S^c) \coloneqq |E_H(S, S^c)|$ and $e_H(S) \coloneqq |E_H(S)|$ and also write $\partial_H(S)=E_H(S, S^c)$ for the edge-boundary of $S$. We denote by $N_H(S)$ the external neighborhood of $S$ in $H$. When the graph $H$ is clear from the context, the subscript may be omitted.
Given a real number $r\in \mathbb{R}$, we denote by $\{r\}$ the fractional part of $r$, which is defined as $\{r\}=r-\lfloor r\rfloor$. For real numbers $x,y,z$, we will write $x = y \pm z$ as a notation for $y-z\le x\le y+z$. Unless otherwise stated, all logarithms are natural logarithms. For $\ell\in \mathbb{N}$, set $[\ell]\coloneqq \{1,2,\ldots,\ell\}$. Throughout the paper, we ignore rounding signs as long as this does not affect the validity of our arguments.

\subsection{Standard Probabilistic Tools}

We will make use of the following fairly standard probability bounds. The first is the well-known Markov's inequality (see, for example, Appendix A in \cite{zbMATH06566409}):
\begin{claim}\label{markov}
    Let $X$ be a non-negative random variable and set $t>0$. Then,
    \begin{align*}
        \mathbb{P}\left(X\ge t\right)\le \frac{\mathbb{E}[X]}{t}.
    \end{align*}
\end{claim}

The second two inequalities are Chernoff-type tail bounds for the binomial and hypergeometric distributions (see, for example, Appendix A in \cite{zbMATH06566409} and \cite[Theorem~2.10]{zbMATH01540669}):
\begin{claim} \label{cher}
    Let $X\sim \mathrm{Bin}(n,p)$. Then for any $0<t\leq \frac{np}{2}$,
    \begin{align*}
        \mathbb{P}(|X-np|\geq t)\leq 2\exp \left(-\frac{t^2}{3np}\right).
    \end{align*}
    Moreover, for any $s>0$, 
    \begin{align*}
        \mathbb{P}(X>s\cdot np)\le \left(\frac{e}{s}\right)^{s\cdot np}.
    \end{align*}
\end{claim}

\begin{claim}\label{cher_hyp}
Fix integers $m,n,N$ and define a random subset $S\subseteq [N]$ distributed uniformly over all subsets of size $n$.
Denote by $X$ the size of $S\cap [m]$.
Then for all $t\ge 0$,
\begin{align*}
\max\{\mathbb{P}(X-\mathbb E X \ge t),\mathbb{P}(X - \mathbb E X \le -t)\} &\le \exp\left(-\frac{t^2}{2\mathbb E X + 2t/3}\right).
\end{align*}
\end{claim}

The last claim is a variant of the well-known Azuma-Hoeffding inequality (see, for example, Chapter 7 in \cite{zbMATH06566409}):
\begin{claim} \label{hoef}
    Let $m \in \mathbb{N}$ and let $p \in [0, 1]$. Let $X = (X_1, X_2, \dots , X_m)$ be a random vector with range $\{0, 1\}^m$ with $X_i$ being independent $Bernoulli(p)$ for every $1\le i\le m$. Let $f:\{0, 1\}^m \xrightarrow{} \mathbb{R}$ be such that there exists a constant $C \in \mathbb{R}$ such that for every $x, x' \in \{0, 1\}^m$ which differ only in one coordinate,
    $$|f(x)-f(x')|\leq C.$$
    Then, for every $t \geq 0$,
    \[
        \mathbb{P}\Big(|f(X)-\mathbb{E}[f(X)]|\geq t\Big)\leq 2\exp \left(-\frac{t^2}{2C^2m}\right).
    \]
\end{claim}

\subsection{Percolation Preliminaries}

We will utilize the following lemma, allowing us to find large matchings in percolated subgraphs; it follows immediately from \cite[Lemma 3.8]{zbMATH07915715}.

\begin{lemma}\label{matching}
Let $H$ be a $d$-regular graph. Let $0<\delta<\frac{1}{2}$ be a constant, and let $s\ge \Omega(d)$. Let $F\subseteq E(H)$ be such that $|F|= s$, and let $q=\frac{\delta}{d}$. Then, $F_{q}$, a random subset of $F$ obtained by retaining each edge independently with probability $q$, contains a matching of size at least $\frac{\delta^2s}{d}$ with probability at least $1-\exp\left(-\frac{\delta^2s}{d}\right)$.
\end{lemma}

We use a slightly strengthened version of a standard elementary claim (see for example Lemma 3.9 in \cite{zbMATH07915715}) and provide a short proof of it for completeness.

\begin{lemma}\label{deg_in_G_p}
    Let $c>0$ be a constant, and let $p=\frac{c}{d}$, where $d = \omega_n(1)$. Let $H$ be a $d$-regular $n$-vertex graph. Then \textbf{whp} $H_p$ contains at most $nd^{-\frac{\log \log d}{4}}$ vertices of degree at least $\log d$.
\end{lemma}

\begin{proof}
    Let $v\in V(H)$. The degree of $v$ in $H_p$ is distributed as $\mathrm{Bin}(d,p)$, hence, by Claim~\ref{cher},
    \begin{align*}
        \mathbb{P}(d_{H_p}(v)\ge \log d)\le \left(\frac{ec}{\log d}\right)^{\log d}\le d^{-\frac{\log \log d}{2}}.
    \end{align*}
    Thus, the expected number of vertices $v\in V(H)$ with $d_{H_p}(v)\ge \log d$ is at most $nd^{-\frac{\log \log d}{2}}$. By using Markov's inequality, we establish the lemma. 
\end{proof}

The following result is due to Frieze, Krivelevich, and Martin \cite[Theorem 1]{zbMATH02037533} and is about the emergence of the giant component of the percolated $(n,d,\lambda)$-graph in the supercritical regime.

\begin{lemma}[{\cite[Theorem 1]{zbMATH02037533}}]\label{FKM}
    Let $\sigma>1$ be a constant. Let $n,d$ be large integers. Assume that $d = \omega_n(1)$ and let $\lambda=o(d)$. Set $p=\frac{\sigma}{d}$. Let $H=(V,E)$ be an $(n,d,\lambda)$-graph. Then, \textbf{whp} $H_p$ contains a unique giant component of asymptotic size $\left(1-\frac{\bar{\sigma}}{\sigma}\right)n$ and the remaining components are of size $O(\log n)$, where $\bar{\sigma}$ is the unique solution (other than $\sigma$) of the equation:
    \[
        xe^{-x}=\sigma e^{-\sigma}.
    \]
\end{lemma}

Next, we prove a quite general result concerning `middle-sized' components in the percolated graph. We will make use of it in the proofs of Theorems \ref{L1} and \ref{setup_2_expansion}. Before stating and proving the lemma, we require a modified version of the classical Breadth First Search (BFS) exploration algorithm (as presented in \cite{zbMATH08098276}).

The algorithm receives as an input a graph $H=(V,E)$ with an ordering $\sigma$ on its vertices, a subset $U\subseteq V$, and a sequence $(X_i)_{i=1}^{|E|}$ of i.i.d Bernoulli$(p)$ random variables. 

We maintain three sets throughout the process: $S$, the set of vertices whose exploration has been completed; $Q$, the set of vertices currently being explored, kept in a queue (first-in-first-out discipline); and $Y$, the set of vertices which have not yet been processed. We initialize $S=\emptyset$, $Q=U$ (with the ordering according to $\sigma$), and $Y=V\setminus U$. The process stops once $Q$ is empty.

At round $t$ of the algorithm execution, we consider the first vertex $v$ in $Q$, and query the first edge between $v$ and $Y$ according to the order $\sigma$ (that is, the edge connecting $v$ to the first -- according to $\sigma$ -- vertex in $Y$). If $X_t=1$, we retain this edge and move its endpoint vertex in $Y$ to $Q$. If $X_t=0$, we discard the edge and continue. If there are no neighbors of $v$ in $Y$ left to query, we move $v$ from $Q$ to $S$ and continue. 

Note that at every moment we have queried (and answered in the negative) all edges between current $S$ and $Y$. Moreover, if $U=\{v\}$ for some $v\in V(H)$, then once $Q$ is empty, $H_p[S]$ has the same distribution as $C_{H_p}(v)$.

\begin{lemma}\label{nointer_general}
    Let $\epsilon,\epsilon',\alpha>0$ be sufficiently small constants satisfying $0<\alpha\le \frac{\epsilon}{4}$ and $\frac{\epsilon}{2}<\epsilon'\le \epsilon$. Also, let $b,\delta>0$ be some constants. Let $n,d$ be large integers. Suppose that $H$ is a $d$-regular $n$-vertex graph satisfying $\mathrm{iso}(b)$ and $\mathrm{EE}(t,\alpha)$ for some $t\coloneqq t(n,d)\in \mathbb{N}$. Set a probability $p=\frac{1+\epsilon'}{d}$, and interval $I=[\delta d,t]$, and a vertex $v\in V(H)$. Then,
    \[
        \mathbb{P}(|C_{H_p}(v)|\in I)\le \exp \left(-\frac{\epsilon^2 \delta}{300}\cdot d\right).
    \]
    Moreover, \textbf{whp} the number of vertices lie in components of size between $\delta d$ and $t$ in $H_p$ is at most $n\cdot \exp \left(-\epsilon^3\delta d\right)$.
\end{lemma}

\begin{proof}
    Let $v\in V(H)$ and $k\in I=\left[\delta d,t\right]$. Utilize the BFS algorithm, with $U=\{v\}$. If $|C_{H_p}(v)|=k$, then there is a round $r$ of the algorithm where $Q$ is empty and $|S|=k$. Since $|S|=k\le t$, and since $H$ satisfies $\mathrm{EE}(t,\alpha)$, at round $r$, we have had at least $(1-\alpha)d|S|$ queries corresponding to $e_H\left(S,S^c\right)$, and at least additional $k-1$ queries corresponding to the internal edges of $C_{H_p}(v)$ which have been explored. Further, only $k-1$ queries were answered in the positive. Therefore,
    \begin{align*}
        \mathbb{P}(|C_{H_p}(v)|=k)\le \mathbb{P}\left(\mathrm{Bin}\left(\left(1-\alpha\right)dk+k-1, p\right)\le k-1\right).
    \end{align*}
    Set $X\sim \mathrm{Bin}\left(\left(1-\alpha\right)dk+k-1, p\right)$.
    Notice that
    \begin{align*}
        k+\frac{\epsilon'-\alpha}{2}\cdot k\le \mathbb{E}[X]=\left(1-\alpha\right)(1+\epsilon')k+\frac{(1+\epsilon')(k-1)}{d}\le k+(\epsilon'-\alpha)k,
    \end{align*}
    where the last inequality holds since $d=\omega(1)$. Hence, by Claim~\ref{cher},
    
    \begin{align*}
        \mathbb{P}(|C_{H_p}(v)|=k)\le \mathbb{P}\left(X\le k\right)\le 2\exp \left(-\frac{(\mathbb{E}[X]-k)^2}{3\mathbb{E}[X]}\right)\le 2\exp \left(-\frac{(\epsilon'-\alpha)^2}{12(1+\epsilon'-\alpha)}\cdot k\right).
    \end{align*}
    Since $\alpha\le \frac{\epsilon}{4}<\frac{\epsilon}{2}<\epsilon'$ and $k\ge \delta d$, we have
    \begin{align*}
        \mathbb{P}(|C_{H_p}(v)|=k)\le 2\exp \left(-\frac{(\epsilon'-\alpha)^2}{12(1+\epsilon'-\alpha)}\cdot k\right)\le \exp\left(-\frac{\epsilon^2}{200}\cdot k\right).
    \end{align*}
    Hence,
    \begin{align*}
        \mathbb{P}\left(|C_{H_p}(v)|\in I\right)&\le \sum_{k=\delta d}^{t} \exp \left(-\frac{\epsilon^2}{200}\cdot k\right)\le \sum_{k=\delta d}^{\infty} \exp \left(-\frac{\epsilon^2}{200}\cdot k\right)\\
        &\le \frac{400}{\epsilon^2}\cdot \exp \left(-\frac{\epsilon^2 \delta}{200}\cdot d\right)\le \exp \left(-\frac{\epsilon^2 \delta}{300}\cdot d\right).
    \end{align*}
    This completes the first part of the lemma.
    Now, for every vertex $v\in V(H)$, set $X_v$ as the indicator of the event that $|C_{H_p}(v)|\in I$, and set $X=\sum_{v\in V(H)} X_v$. Since
    \begin{align*}
        \mathbb{E}[X]\le n\cdot \exp \left(-\frac{\epsilon^2 \delta}{300}\cdot d\right),
    \end{align*}
    then by Claim~\ref{markov}, and since $d=\omega(1)$,
    \begin{align*}
        \mathbb{P}\left(X\ge n\cdot \exp \left(-\epsilon^3\delta d\right)\right)\le \exp\left(\epsilon^3\delta d-\frac{\epsilon^2 \delta}{300}\cdot d\right)=o(1),
    \end{align*}
    which establishes the lemma.
\end{proof}

Finally, recall that we denote by $\gamma(\epsilon)$ the survival probability of the $\mathrm{Po}(1+\epsilon)$ branching process, where $\epsilon>0$ is a small constant. Let $\delta>0$ be a sufficiently small constant in terms of $\epsilon$. In the proofs of Theorems \ref{L1} and \ref{general} we use the fact that
\begin{align}\label{gamma}
    \gamma(\epsilon)-\gamma(\epsilon-\delta)\le 2\delta.
\end{align}

\subsection{Multiple-Exposure/Sprinkling Argument}\label{sprinkling sec}

A (by now) standard tool in percolation theory is a multiple-exposure/sprinkling argument \'{a} la Ajtai-Koml\'{o}s-Szemer\'{e}di \cite{zbMATH03769676}. In the general case, given a $d$-regular $n$-vertex graph $H$, a positive integer $k$ and a probability $p\in [0,1]$, we define $p_1, p_2,\ldots, p_k\in [0,1]$ such that
\begin{align}\label{sprinkling}
    1-p=(1-p_1)(1-p_2)\cdot \ldots \cdot (1-p_k).
\end{align} 
Now, in addition to the graph $H_p$, we consider the graphs $H_{p_1},H_{p_2},\ldots, H_{p_k}$. Furthermore, inductively, we set $H_1=H_{p_1}$ and for every $1<j\in [k]$, we set $H_j= H_{j-1}\cup H_{p_j}$. Observe that by \eqref{sprinkling} $H_p$ has the same distribution as $H_k=H_{p_1}\cup H_{p_2}\cup \ldots \cup H_{p_k}$. Usually, we set $k=2$ or $k=3$, take $p_1$ to be close to $p$ and $p_2,\ldots,p_k$ relatively small, and think of $H_{p_2},\ldots, H_{p_k}$ as graphs whose role is to fix $H_{p_1}$ in order to reach our goal on the graph $H_p$. We use the tool of multiple-exposure/sprinkling in all our results in Section \ref{main results}, where we consider $p=\frac{1+\epsilon}{d}$, where $\epsilon>0$ is a sufficiently small constant. 

In this section, we set $k=3$ and $p=\frac{1+\epsilon}{d}$, where $\epsilon>0$ is a sufficiently small constant. We further define $p_2=p_3=\frac{\delta}{d}$, where $\delta>0$ is a sufficiently small constant in terms of $\epsilon>0$ and $p_1$ is taken such that $1-p=(1-p_1)(1-p_2)(1-p_3)$. Let us note that the following results remain valid for $k=2$ and $p_2=\frac{2\delta}{d}$.

Below we state auxiliary results from \cite{zbMATH08225676} that are valid for any $d$-regular $n$-vertex graph $H$, assuming that $d=\omega_n(1)$ (without any further assumptions). We begin with an exponentially unlikely event in $H_1$.

\begin{lemma}[{\cite[Lemma 4.4]{zbMATH08225676}}]\label{unlikely_event_H1}
    Let $v\in V(H)$. The probability that there exists $A\subseteq N_H(v)$ such that
    \begin{align*}
        \left|\bigcup_{u\in A}C_{H_1}(u)\right|\in \left[\delta d, 2\delta d\right]
    \end{align*}
    and $|A|\le \delta^2d$ is at most $\exp\left(-\delta^2 d\right)$.
\end{lemma}

Before stating the next two results, let us introduce two subsets of vertices of $V(H)$. Let $j\in \{1,2,3\}$. The first set contains vertices that lie in `large' components of $H_j$:
\begin{align*}
    V_L(H_j):=\left\{v\in V(H)\;|\; |C_{H_j}(v)|\ge \delta d\right\}.
\end{align*}
The second set consists of vertices that are adjacent to many vertices in `large' components of $H_j$:
\begin{align*}
    W_L(H_j):=\left\{w\in V(H)\;|\; \left|N_H(w)\cap V_L(H_j)\right|\ge \delta^5 d\right\}.
\end{align*}

The following two results determine the number of vertices in $H_j$ ($j=1,2,3$) that lie in `large' components (Lemma~\ref{volume}), and the number of vertices in $H$ that are adjacent to many vertices in `large' components (Lemma~\ref{W}). Recall that $\gamma(\epsilon)$ the survival probability of the $\mathrm{Po}(1+\epsilon)$ branching process.

\begin{lemma}[{\cite[Lemma 4.3]{zbMATH08225676}}]\label{volume}
    For every $j\in \{1,2,3\}$, \textbf{whp}
    \[
        |V_L(H_j)|\in \left[(1-o_d(1))\gamma(\epsilon-2\delta)n,(1+o_d(1))\gamma(\epsilon)n\right].
    \]
\end{lemma}

\begin{lemma}[{\cite[Lemma 4.5]{zbMATH08225676}}]\label{W}
    \textbf{Whp} $|W_L(H_2)|\ge (1-\exp(-\delta^6d))n$.
\end{lemma}

\subsection{Implications of Good Vertex Expansion}

In this section we show that good vertex expansion guarantees the existence of certain combinatorial substructures in the graph. The first result show that under some vertex expansion assumptions, a long cycle appears.

\begin{claim}[{\cite[Theorem 1]{zbMATH07058289}}]\label{cycle}
    Let $k\ge 1,\ell\ge 2$ be integers. Let $H$ be a graph on more than $k$ vertices such that for every $U\subseteq V(H)$ with $\frac{k}{2}\le |U|\le k$, $|N_H(U)|\ge \ell$. Then $H$ contains a cycle of length at least $\ell+1$.
\end{claim}

The next claim allows us to deduce the existence of a large complete minor in a graph without any small separators. Since graphs with good vertex expansion has no small separator (see, for example, \cite[Section 5]{zbMATH07307104}) the following claim essentially links vertex expansion and the existence of a large complete minor. 

\begin{claim}[{\cite[Theorem 1.2]{zbMATH08033971}}]\label{minor}
    Let $H$ be a graph with $n$ vertices and with no $K_t$-minor. Then there is $S\subseteq V(H)$ with $|S|=O\left(t\sqrt{n}\right)$ such that each connected component of $H\setminus S$ has at most $\frac{2}{3}n$ vertices.
\end{claim}

\subsection{Further Well-Known Results}

Throughout this paper we also use two more classic and well-known results. The first one is an upper bound on the number of $k$-vertex trees in a graph $H$, which follows immediately from \cite[Lemma 2]{zbMATH01288356}:
\begin{claim}\label{trees}
    Let $H$ be a graph of maximum degree at most $d$, let $v \in V (H)$, and let $k \in \mathbb{N}$. Let $t_k(H, v)$ denote the number of trees on $k$ vertices in $H$ rooted at $v$. Then
    \begin{align*}
        t_k(H, v)\le \frac{k^{k-2}d^{k-1}}{(k-1)!}\le (ed)^{k-1}.
    \end{align*}
\end{claim}

The second result is a well-known fundamental result by Alon and Milman \cite{zbMATH03875335}. This result concerns pseudo-random graphs which will be a cornerstone in our constructions presented in the proofs of Theorems \ref{no_L1}, \ref{no_expansion}, and \ref{general_large_diam}. We give here a version of this result that will be useful for us:

\begin{claim}[{\cite{zbMATH03875335}}]\label{eml}
    Let $H = (V, E)$ be a $d$-regular $n$-vertex graph, and denote by $\lambda_2$ the second eigenvalue of its adjacency matrix. Then for every $U\subseteq V$,
    \begin{align*}
        e_H(U,U^c)\ge \frac{d-\lambda_2}{n}|U|(n-|U|).
    \end{align*}
\end{claim}

\section{Guaranteeing the Giant Component --- Proof of Theorem~\ref{L1}}\label{pL1}

Throughout this section we use notations and assumptions in Theorem~\ref{L1}. We further assume that $d=o\left(\frac{n}{\log n}\right)$. The complementary assumption requires essentially the same treatment, with the necessary modifications. See Section 4 in \cite{zbMATH08242616} for further details on this case. Set $\delta=\frac{1}{(Cb)^{1/20}}$. Notice that since $Cb$ is sufficiently large in terms of $\epsilon$, we may assume that $\delta\le \epsilon^3$. We utilize a multiple-exposure/sprinkling argument \'{a} la Ajtai-Koml\'{o}s-Szemer\'{e}di \cite{zbMATH03769676}. Let $p_2=p_3=\frac{\delta}{d}$, and let $p_1$ be such that $(1-p_1)(1-p_2)(1-p_3)=1-p$. We define the following graphs:

\begin{itemize}
    \item $G_1=G_{p_1}$. Set $\rho_1=p_1$.
    \item $G_2=G_1\cup G_{p_2}$. Set $\rho_2$ such that $(1-p_1)(1-p_2)=(1-\rho_2)$. Note that $G_2=G_{\rho_2}$.
    \item $G_3=G_2\cup G_{p_3}$. Set $\rho_3=p$. Note that $G_3=G_{\rho_3}$.
\end{itemize}

Note that $G_p$ has the same distribution as $G_3=G_{p_1} \cup G_{p_2} \cup G_{p_3}$. Further notice that for every $j\in \{1,2,3\}$, there is $\epsilon-2\delta\le\epsilon_j\le \epsilon$, such that $\rho_j=\frac{1+\epsilon_j}{d}$. For sufficiently small $\delta$, we can assume that $\frac{\epsilon}{2}<\epsilon_j\le \epsilon$.

Following Section \ref{sprinkling sec}, for each $j\in \{1,2,3\}$, we consider two subsets of vertices. The first one contains vertices that lie in `large' components of $G_j$:
\begin{align*}
    V_L(G_j)=\left\{v\in V(G)\;|\; |C_{G_j}(v)|\ge \delta d\right\}.
\end{align*}
The second set consists of vertices that are adjacent to many vertices in `large' components of $G_j$:
\begin{align*}
    W_L(G_j)=\left\{w\in V(G)\;|\; \left|N_G(w)\cap V_L(G_j)\right|\ge \delta^5 d\right\}.
\end{align*}

In the proof of Theorem~\ref{L1}, we rely on a couple of lemmas from \cite{zbMATH08225676} that determine the typical number of vertices in $G_j$ ($j=1,2,3$) that lie in `large' components (Lemma~\ref{volume}) and the number of vertices in $G$ that are adjacent to many vertices in `large' components (Lemma~\ref{W}). Our analysis counts the number of vertices that lie in `middle-sized' components (Lemma~\ref{nointer}), which allows us to carefully control a union bound argument presented at the end of the section.

Since for every $1\le j\le 3$, we have $\frac{\epsilon}{2}<\epsilon_j\le \epsilon$, and since our graph $G$ satisfies $\mathrm{iso}(b)$ and $\mathrm{EE}(Cd,\alpha)$, we can apply Lemma~\ref{nointer_general} which we have proven in Section \ref{pre} to obtain:

\begin{lemma}\label{nointer}
    For every $j\in \{1,2,3\}$, \textbf{whp} the number of vertices in components of size between $\delta d$ and $Cd$ in $G_j$ is at most $n\cdot \exp \left(-\epsilon^3\delta d\right)$.
\end{lemma}

Now, we show that after sprinkling with probability $p_3$, \textbf{whp} all but at most $\delta n$ of the vertices in $V_L(G_2)$ merge into a single component, using Lemmas \ref{matching}, \ref{volume}, \ref{W}, and \ref{nointer}.

\begin{proof}[Proof of Theorem~\ref{L1}]
    By Lemma~\ref{volume} and \eqref{gamma}, we have
    \begin{align*}
        |V_L(G_3)|\le (1+o_d(1))\gamma(\epsilon)n, \quad \quad |V_L(G_2)|\ge (1-o_d(1))\gamma(\epsilon-2\delta)n\ge \gamma(\epsilon)n-5\delta n.
    \end{align*}
    As mentioned above, all that remains is to show that after sprinkling with probability $p_3$, \textbf{whp} all but at most $\delta n$ of the vertices in $V_L(G_2)$ merge into a single component. Let $A \sqcup B$ be a partition of the components forming $V_L(G_2)$, satisfying $\delta n\le k\coloneqq |A|\le |B|$. Let $A'$ be $A$ together with all the vertices in $V(G)\backslash B$ which have at least $\frac{\delta^5 d}{2}$ neighbors in $A$, and let $B'$ be $B$ together with all the vertices in $V(G)\backslash A'$ which have at least $\frac{\delta^5 d}{2}$ neighbors in $B$. Note that $A'\cap B'=\emptyset$ and $V_L(G_2)\cup W_L(G_2)\subseteq A'\cup B'$. Furthermore, $|A'|,|B'|\ge k$ and one of them is of order at most $\frac{n}{2}$. We may assume without loss of generality that $k\le |A'|\le \frac{n}{2}$. Then, since $G$ satisfies $\mathrm{iso}(b)$, we have that $e_G(A',(A')^c)\ge b|A'|\ge bk$. Also, by Lemma~\ref{W}, \textbf{whp} $|V(G)\backslash W_L(G_2)|\le \exp(-\delta^6d)n=o(n/d)$, therefore $e_G(A',B')\ge bk/2$.
    
    Let us now expose these edges (between $A'$ and $B'$) with probability $p_3$. 
    By Lemma~\ref{matching}, with probability at least $1-\exp \left(-\frac{\delta^2bk}{2d}\right)$ there is a matching of size at least $\frac{\delta^2bk}{2d}$ between $A'$ and $B'$ in $G_{p_3}$. 
    By the definition of $A'$ and $B'$, from the endpoints of each edge in the matching there are at least $\frac{\delta^5 d}{2}$ edges into $A$, and at least $\frac{\delta^5 d}{2}$ edges into $B$. 
    The probability that a given edge in the matching does not extend to a path of length three between $A$ and $B$ is thus at most
    \begin{align*}
        1-\left(1-(1-p_3)^{\frac{\delta^5 d}{2}}\right)^2&= \left(2-(1-p_3)^{\frac{\delta^5 d}{2}}\right)(1-p_3)^{\frac{\delta^5 d}{2}}\le \left(2-\left(1-\frac{\delta^5 d}{2}\cdot \frac{\delta}{d}\right)\right)\exp\left(-\frac{\delta^6}{2}\right)\\
        &= \left(1+\frac{\delta^6}{2}\right)\exp\left(-\frac{\delta^6}{2}\right)\le \exp\left(-\delta^{13}\right).
    \end{align*}
    These events are independent for the edges of the matching, and thus
    \begin{align}\label{eq: no path A,B}
        \mathbb{P}(\text{there is no path in $G_{p_3}$ between $A$ and $B$})\le \exp\left(-\frac{\delta^{13}\cdot \delta^2bk}{2d}\right)=\exp\left(-\frac{\delta^{15} bk}{2d}\right).
    \end{align}

    We now bound from above the number of partitions of the components forming $V_L(G_2)$ into $V_L(G_2)= A \sqcup B$ with $|A| = k$. We can write $V_L(G_2)=M\sqcup L$, where
    \begin{align*}
        M=\left\{v\in V(G)\;|\; |C_{G_2}(v)|\in [\delta d,Cd]\right\}\quad \text{and}\quad L=\left\{v\in V(G)\;|\; |C_{G_2}(v)|> C d\right\}.
    \end{align*}
    By Lemma~\ref{nointer}, \textbf{whp} $|M|\le n\cdot \exp\left(-\epsilon^3\delta d\right)$. For every partition of the components forming $V_L(G_2)$ into $V_L(G_2)= A \sqcup B$, we have to divide both $M$ and $L$ into $A \sqcup B$. The number of partitions of $M$ into two parts is at most 
    \begin{align*}
        2^{|M|}\le \exp\left(n\cdot \exp\left(-\epsilon^3\delta d\right)\right).
    \end{align*}
    Hence, the number of ways to choose $A$ (and hence the partition) is at most 
    \begin{align}\label{eq: partition}
        2^{|M|}\cdot \sum_{i=1}^{k/(Cd)}\binom{n/(Cd)}{i}\le \exp\left(n\cdot \exp\left(-\epsilon^3\delta d\right)\right)\cdot \left(\frac{en}{k}\right)^\frac{k}{Cd},
    \end{align}
    since every vertex $v\in L$ is part of a component of size greater than $Cd$ in $G_2$. Thus, by the union bound and using \eqref{eq: no path A,B} and \eqref{eq: partition}, the probability that there exists such a partition $\{A,B\}$ without a path in $G_{p_3}$ between $A$ and $B$ is at most
    \begin{align}\label{eq: final_prob}
        \sum_{k=\delta n}^{n/2} \exp\left(n\cdot \exp\left(-\epsilon^3\delta d\right)\right)&\cdot \left(\frac{en}{k}\right)^\frac{k}{Cd}\cdot \exp\left(-\frac{\delta^{15} bk}{2d}\right).
    \end{align}
    Observe that
    \begin{align*}
        \sum_{k=\delta n}^{n/2} \left(\frac{en}{k}\right)^\frac{k}{Cd}\cdot \exp\left(-\frac{\delta^{15} bk}{2d}\right)&=\sum_{k=\delta n}^{n/2} \exp\left(\frac{k}{Cd}\log \left(\frac{en}{k}\right) -\frac{\delta^{15} bk}{2d}\right)\\
        &\le \sum_{k=\delta n}^{n/2} \exp\left(\frac{k}{d}\left(\frac{1+\log \frac{1}{\delta}}{C}-\frac{\delta^{15}b}{2}\right)\right),
    \end{align*}
    where we used $k\ge \delta n$. Now, by our choice of $\delta=\frac{1}{(Cb)^{1/20}}$, the first summand is the largest term on the right-hand side of the above inequality. Hence
    \begin{align*}
        \sum_{k=\delta n}^{n/2} \exp\left(\frac{k}{d}\left(\frac{1+\log \frac{1}{\delta}}{C}-\frac{\delta^{15}b}{2}\right)\right)\le n\exp\left(\frac{\delta n}{d}\left(\frac{1+\log \frac{1}{\delta}}{C}-\frac{\delta^{15} b}{2}\right)\right).
    \end{align*}
    Altogether, \eqref{eq: final_prob} can be bounded from above by
    \begin{align*}
        \sum_{k=\delta n}^{n/2} \exp\left(n\cdot \exp\left(-\epsilon^3\delta d\right)\right)&\cdot \left(\frac{en}{k}\right)^\frac{k}{Cd}\cdot \exp\left(-\frac{\delta^{15} bk}{2d}\right)\\
        &\le \exp\left(n\cdot \exp\left(-\epsilon^3\delta d\right)\right)\cdot n\exp\left(\frac{\delta n}{d}\left(\frac{1+\log \frac{1}{\delta}}{C}-\frac{\delta^{15} b}{2}\right)\right)=o(1),
    \end{align*}
    where we used $d=o\left(\frac{n}{\log n}\right)$. We conclude that \textbf{whp} $G_p$ contains a component of size at least $\gamma(\epsilon)n-5\delta n-\delta n\ge \left(\gamma(\epsilon)-\frac{10}{(Cb)^{1/20}}\right)n$. Furthermore, with the same probability, the second-largest component is of size at most $(1+o_d(1))\frac{10}{(Cb)^{1/20}}n$.
\end{proof}

Let us note that most of the constants in the proof (such as 1/20 and the exponent 15 in $\delta^{15}$) are not optimal, and can be improved with a more tedious implementation.

\section{Vertex Expansion of Linear-Sized Sets --- Proof of Theorem~\ref{general}}\label{Pgeneral}

Throughout this section we use notations and assumptions in Theorem~\ref{general}. We further assume that $d=o\left(\frac{n}{\log n}\right)$, as otherwise the conclusion of Theorem~\ref{general} becomes trivial.

We utilize a double-exposure/sprinkling argument. Let $\delta\coloneqq \delta(\epsilon)\le \epsilon^3$. Define $p_2=\frac{\delta}{d}$ and let $p_1$ be such that $(1-p_1)(1-p_2)=1-p$, noting that $p_1\geq \frac{1+\epsilon-\delta}{d}$. Notice that $G_p$ has the same distribution as $G_{p_1}\cup G_{p_2}$. We abbreviate $G_1\coloneqq G_{p_1}$ and $G_2\coloneqq G_{p_1}\cup G_{p_2}$. Denote by $L_1'$ and $L_1$ the largest components of $G_{p_1}$ and $G_p$, respectively. Furthermore, if we let $\gamma$ and $\gamma_1$ be the survival probabilities of the $\mathrm{Po}(1+\epsilon)$ and $\mathrm{Po}(dp_1)$ branching processes, respectively, and since $dp_1\ge 1+\epsilon-\delta$, then we may assume by \eqref{gamma} that
\begin{align*}
    \gamma-\gamma_1\le 3\delta.
\end{align*}
Furthermore, by Theorem~\ref{ERCP1} (with Remark \ref{re: sprinkling_ercp1}), \textbf{whp}
\begin{align}\label{eq: pre_big_components}
    |V(L_1)|=(\gamma\pm \epsilon^3)n,\quad |V(L_1')|=(\gamma_1\pm (\epsilon-\delta)^3)n.
\end{align}
Also observe that \textbf{whp} $L_1'\subseteq L_1$ since \textbf{whp} in $G_2$, the size of the second largest component is smaller than $|V(L_1')|$. Combining \eqref{eq: pre_big_components} with $\delta\le \epsilon^3$ yields
\begin{align}\label{big_components}
    \big||V(L_1)|-|V(L_1')|\big|\le 5\epsilon^3 n.
\end{align}

The proof proceeds as follows: we start with large sets inside $L_1'$. Lemma~\ref{key-lemma} shows that \textbf{whp} such sets expand themselves well in one of two ways: either they have good vertex expansion within $L_1'$, or many vertex-disjoint paths leave them in $G_{p_2}$. The second type of expansion is guaranteed due to Lemma~\ref{paths_setup2}. Finally, we look at large sets in $L_1$. We finish due to \eqref{big_components} and since we deal with large enough sets. 

The following argument shows that with very high probability one can find many vertex-disjoint paths between the two parts of every partition of $V(L_1')$ into two parts that are not too small. It follows from the proof of Theorem 1 in \cite{zbMATH08242616} (particularly from Lemma 3.5 there):

\begin{lemma}\label{paths_setup2}
    There exist constants $C'\coloneqq C'(\epsilon)>0$ and $c = c(\epsilon, c_1,c_2,c_3,\alpha)>0$ such that the following holds. Let $A\cup B=V(L_1')$ be a partition of $V(L_1')$ with $\min\{|A|,|B|\}=k$ where $k\ge C'\log n$. Then, with probability at least $1-\exp\left(-\frac{ck}{d}\right)$, there exists a family of $\frac{ck}{d}$ many vertex-disjoint $A-B$-paths of length at most $2(\log_d\log n+1)+1$ in $G_{p_2}$.
\end{lemma}

In the following lemma, we utilize Lemma~\ref{paths_setup2} to show that any typical large subset of $L_1'$ expands itself well.

\begin{lemma}\label{key-lemma}
    There exists a constant $c>0$ such that \textbf{whp} for any subset of the largest component $S\subseteq V(L_1')$ with $\frac{2\epsilon^2}{3} n\le |S|\le \frac{3\epsilon n}{2}$, either
    \begin{align*}
        |N_{L_1'}(S)|\ge \frac{c}{d\log d}|S|,
    \end{align*}
    or there exists a family of at least $\frac{c|S|}{d}$ many vertex-disjoint paths of length at most $2(\log_d\log n+1)+1$ between $S$ and $V(L_1')\setminus S$ in $G_{p_2}$.
\end{lemma}
We note that the assumption that $|S| \ge \frac{2\epsilon^2}{3} n$ can be relaxed, however it suffices for our usage and allows for a simpler proof. The proof of Lemma~\ref{key-lemma} is inspired by the proof of Lemma 5.3 in \cite{zbMATH07915715}.

\begin{proof}

By Lemma~\ref{deg_in_G_p}, \textbf{whp} there are at most $\frac{n}{d^4}$ vertices with degree larger than $\log d$. We continue assuming that this holds deterministically.

Given a set $S\subseteq V(L_1')$ with $\frac{2\epsilon^2}{3} n\le |S|\le \frac{3\epsilon n}{2}$, let $k\coloneqq |S|$ and let $b\coloneqq |N_{L_1'}(S)|$. As we aim to bound the vertex expansion of the set $S$, we may assume that $b< \frac{ck}{d\log d}$, as otherwise the claim holds. In order to facilitate a union bound argument, let us estimate the number of subsets $S$ of size $k$ in $L_1'$ such that $|N_{L_1'}(S)|=b$. 

Let $e = |\partial_{L_1'}(N_{L_1'}(S))|$. Since there are at most $\frac{n}{d^4}$ vertices with degree larger than $\log d$, we have that $e\le \frac{n}{d^3}+b_1\log d$. Since $L_1'$ is connected, there are at most $e$ components in $L_1'\setminus N_{L_1'}(S)$. Furthermore, since $S$ has no neighbors outside $N_{L_1'}(S)$, it must be the union of components in $L_1'\setminus N_{L_1'}(S)$. Hence, the number of ways to choose such an $S$ is at most
\[
    \sum_{e=1}^{\frac{n}{d^3}+b_1\log d} \binom{n}{b}\cdot 2^e\le \binom{n}{b}\cdot 2^{\frac{n}{d^3}+b_1\log d+1}.
\]
Therefore, the number of sets $S\subseteq V(L_1')$ with $|N_{L_1'}(S)|<\frac{ck}{d\log d}$ is at most
\begin{align*}
    \sum_{b=1}^{\frac{ck}{d\log d}}\binom{n}{b}2^{\frac{n}{d^3}+b\log d+1}&\le \left(\frac{en}{\frac{ck}{d\log d}}\right)^{\frac{ck}{d\log d}}2^{\frac{2ck}{d}}\\
    &\le \exp\left(\frac{ck}{d\log d}\left(\log\left(\frac{end\log d}{ck}\right)+2\log d\right)\right)\\
    &\le\exp\left(\frac{2ck}{d\log d}\left(\log\left(\frac{n}{k}\right)+2\log d\right)\right)\le \exp\left(\frac{6ck}{d}\right),
\end{align*}
where we used the fact that $k\ge \frac{2\epsilon^2}{3} n$ in the first and last inequalities.

Fix $S\subseteq V(L_1')$ with $|S|=k$. By Lemma~\ref{paths_setup2}, there is a constant $c'>0$ such that with probability at least $1-\exp\left(-\frac{c'k}{d}\right)$, there exists a family of at least $\frac{c'k}{d}$ vertex-disjoint paths of length at most $2(\log_d\log n+1)+1$ between $S$ and $V(L_1')\setminus S$ in $G_{p_2}$.

The probability there is a set $S$ of size $k$ violating the statement is then at most
\begin{align*}
    \exp\left(\frac{6ck}{d}-\frac{c'k}{d}\right)=\exp\left(\frac{k}{d}\left(6c-c'\right)\right).
\end{align*}

Recall that $k=|S|\ge \frac{2\epsilon^2}{3} n$, and that $d=o\left(\frac{n}{\log n}\right)$. Hence, by choosing $c$ sufficiently small with respect to $c'$, by the union bound the probability of having a set $S$ violating the statement is $o(1)$.
\end{proof}

\begin{proof}[Proof of Theorem~\ref{general}]
    Let $S\subseteq V(L_1)$ such that $\epsilon ^2 n\le |S|\le \frac{3\epsilon n}{2}$. Recalling that \textbf{whp} $L_1'\subseteq L_1$, we define $S_1=S\cap V(L_1')$ and $S_2=S\cap \left(V(L_1)\setminus V(L_1')\right)$. Recall that by \eqref{big_components}, $\big||V(L_1)|-|V(L_1')|\big|\le 5\epsilon^3n$, and so it follows that $|S_1|\ge \frac{2|S|}{3}$. Thus, by Lemma~\ref{key-lemma}, there is a constant $c>0$ such that \textbf{whp} either
    \begin{align*}
        |N_{L_1'}(S_1)|\ge \frac{c|S_1|}{d\log d}\ge \frac{\beta}{d\log d}|S|,
    \end{align*}
    or there is a family of at least $\frac{c|S_1|}{d}\ge\frac{\beta|S|}{d}$ vertex-disjoint paths from $S_1$ to $V(L_1')\setminus S_1 \subseteq V(L_1) \setminus S$. However, since each such path contributes a unique vertex to the neighborhood of $S$ in $V(L_1)$ (the first vertex along the path which is not in $S$), in the latter case $|N_{G_p}(S)|\ge \frac{\beta|S|}{d}$, and so the result follows. For the first case, since $N_{L_1'}(S_1)\subseteq N_{G_p}(S)$, then again, the result follows.
\end{proof}

\begin{proof}[Proof of Corollary~\ref{full expansion in subgraph}]
The proof follows the lines of that of Lemma 2.4 in \cite{zbMATH07307104} or that of Theorem 5 in \cite{zbMATH07915715}.
\end{proof}

\begin{proof}[Proof of Corollary~\ref{properties_general}]
    
    \textbf{Long cycles:} Applying Theorem~\ref{general}, we deduce from Claim~\ref{cycle} that $L_1$ contains a cycle of length $\Omega\left(\frac{n}{d\log d}\right)$.

    \textbf{Hadwiger number:} Let $t\in \mathbb{N}$ and assume that $L_1$ contains no $K_t$-minor. By Claim~\ref{minor}, $L_1$ contains a subset $X\subseteq V(L_1)$ of size $O\left(t\sqrt{|V(L_1)|}\right)\le O\left(t\sqrt{n}\right)$ such that all components of $L_1\setminus X$ are of size at most $\frac{2}{3}|V(L_1)|$.

    It follows that there is some subset $S \subseteq V(L_1)$, which is the union of some components of $L_1\setminus X$, such that $\frac{|V(L_1)|}{3}\le |S|\le \frac{|V(L_1)|}{2}$ and $N_{G_p}(S) \subseteq X$. Hence, by Theorem~\ref{general}
  
    \begin{align*}
        \frac{\beta}{d\log d}|S|\le |N_{G_p}(S)|\le |X|=O\left(t\sqrt{n}\right).
    \end{align*}
    Therefore, since $|S|\ge \frac{1}{6}|V(L_1)|=\Omega(n)$, we have $t=\Omega\left(\frac{\sqrt{n}}{d\log d}\right)$.
\end{proof}

\section{Vertex Expansion of Connected Sets --- Proof of Theorem~\ref{setup_2_expansion}} \label{Psetup_2_expansion}

Throughout this section we use notations and assumptions of Theorem~\ref{setup_2_expansion}. We utilize a multiple-exposure/sprinkling argument. Define $p_2=p_3=\frac{\epsilon^3}{d}$ and let $p_1$ be such that $(1-p_1)(1-p_2)(1-p_3)=1-p$, noting that $p_1\geq \frac{1+\epsilon-2\epsilon^3}{d}$. Also, let $p'$ such that $(1-p')(1-p_3)=1-p$. Notice that $G_{p'}$ has the same distribution as $G_{p_1}\cup G_{p_2}$ and that $G_p$ has the same distribution as $G_{p_1}\cup G_{p_2}\cup G_{p_3}$. We abbreviate $G_1\coloneqq G_{p_1}$, $G_2\coloneqq G_{p_1}\cup G_{p_2}$ and $G_3\coloneqq G_{p_1}\cup G_{p_2}\cup G_{p_3}$. Denote by $L_1'$ and $L_1$ the largest components of $G_2$ and $G_3$, respectively. Note that \textbf{whp} the size of the second largest component in $G_3$ is strictly smaller than $|V(L_1')|$, and hence \textbf{whp} $L_1'\subseteq L_1$.  

Let us discuss the proof outline. The proof starts in a manner similar to the proof of Theorem~\ref{general}. We start with connected sets inside $L_1'$. Lemma~\ref{key lemma} shows that \textbf{whp} such sets expand themselves well in one of two ways: either they have good vertex expansion within $L_1'$, or many vertex-disjoint paths leave them in $G_{p_3}$. The second type of expansion is guaranteed due to Lemma~\ref{paths_setup2}. Now, we look at connected sets in $L_1$. When there is enough mass of our set in $L_1'$ we finish by Lemma~\ref{key lemma} (similarly to the proof of Theorem~\ref{general}). When this is not the case, we show (Lemma~\ref{key lemma 2}) that \textbf{whp} there is a good vertex expansion in $L_1$. This is done due to a matching in $G_{p_3}$ between a subset of $V(L_1')$ and a subset of $V(L_1)\setminus V(L_1')$, which is found \textbf{whp} in Lemma~\ref{middle sets}.

\begin{lemma}\label{key lemma}
For any constant $t>0$, there exist constants $\hat{c}=\hat{c}(C_1,c_1,c_2,c_3,\alpha,\epsilon)>0$ and $C_2\coloneqq C_2(C_1,c_1,c_2,c_3,\alpha,\epsilon)>0$ such that \textbf{whp} for any subset of the largest component $S\subseteq V(L_1')$ with $t\cdot nd^{-C_1}\le |S|\le \frac{3\epsilon n}{2}$ and $|S|\ge C_2 (d\log d+\log n)$, then either
    \begin{align*}
        |N_{L_1'}(S)|\ge \frac{\hat{c}}{d\log d}|S|,   
    \end{align*}    
    or there exists a family of at least $\frac{\hat{c}|S|}{d}$ many vertex-disjoint paths of length at most $2(\log_d\log n+1)+1$ between $S$ and $V(L_1')\setminus S$ in $G_{p_3}$.
\end{lemma}

We prove Lemma~\ref{key lemma} in a manner similar to the proof of Lemma~\ref{key-lemma}.

\begin{proof}
By Lemma~\ref{deg_in_G_p}, \textbf{whp} there are at most $nd^{-\frac{\log \log d}{4}}$ vertices with degree larger than $\log d$. We continue assuming that this holds deterministically.

Given $t\cdot nd^{-C_1}\le |S|\le \frac{3\epsilon n}{2}$ and $|S|\ge C_2 (d\log d+\log n)$, let $k\coloneqq |S|$ and let $b_1\coloneqq |N_{L_1'}(S)|$. As we aim to bound the vertex expansion of the set $S$ in $L_1\supsetneq L_1'$, we may assume that $b_1< \frac{\hat{c}k}{d\log d}$, as otherwise the claim holds. In order to facilitate a union bound argument, let us estimate the number of subsets $S$ of size $k$ in $L_1'$ such that $|N_{L_1'}(S)|=b_1$. Let $e_1 = |\partial_{L_1'}(N_{L_1'}(S))|$. Since there are at most $nd^{-\frac{\log \log d}{4}}\le t\cdot nd^{-C_1-3}$ vertices with degree larger than $\log d$, we have that $e_1\le t\cdot nd^{-C_1-2}+b_1\log d$. Since $L_1'$ is connected, there are at most $e_1+1$ components in $L_1'\setminus N_{L_1'}(S)$. Furthermore, since $S$ has no neighbors outside $N_{L_1'}(S)$, it must be the union of components in $L_1'\setminus N_{L_1'}(S)$. Hence, the number of ways to choose such an $S$ is at most $\binom{n}{b_1}\cdot 2^{e_1+1}$. Thus, there are at most
\begin{align*}
    \sum_{b_1=1}^{\frac{\hat{c}k}{d\log d}}\binom{n}{b_1}2^{t\cdot nd^{-C_1-2}+b_1\log d+1}&\le \left(\frac{en}{\frac{\hat{c}k}{d\log d}}\right)^{\frac{\hat{c}k}{d\log d}}2^{\frac{2\hat{c}k}{d}}\\
    &\le \exp\left(\frac{\hat{c}k}{d\log d}\left(\log\left(\frac{end\log d}{\hat{c}k}\right)+2\log d\right)\right)\\
    &\le\exp\left(\frac{2\hat{c}k}{d\log d}\left(\log\left(\frac{n}{k}\right)+2\log d\right)\right)\le \exp\left(\frac{2(C_1+3)\hat{c}k}{d}\right)
\end{align*}
sets $S\subseteq V(L_1')$ with $|N_{L_1'}(S)|<\frac{\hat{c}k}{d\log d}$, where we used the fact that $k\ge t\cdot nd^{-C_1}$ in the first and last inequalities.

Fix $S\subseteq V(L_1')$ with $|S|=k$. Since $G$ satisfies the assumptions of Theorem~\ref{ERCP1}, we can apply Lemma~\ref{paths_setup2}. Hence, with probability at least $1-\exp\left(-\frac{ck}{d}\right)$, there exists a family of at least $\frac{ck}{d}$ vertex-disjoint paths of length at most $2(\log_d\log n+1)+1$ between $S$ and $V(L_1')\setminus S$ in $G_{p_3}$, where $c$ is the constant from Lemma~\ref{paths_setup2}.

The probability there is a set $S$ of size $k$ violating the statement is then at most
\begin{align*}
    \exp\left(\frac{2(C_1+3)\hat{c}k}{d}-\frac{ck}{d}\right)=\exp\left(\frac{k}{d}\left((2C_1+6)\hat{c}-c\right)\right).
\end{align*}

Recall that $k=|S|\ge C_2 (d\log d+\log n)$ and that $d=\omega(1)$. Hence, for $\hat{c}$ small enough with respect to $c$ and $C_1$, we get $c'\coloneqq (2C_1+6)\hat{c}-c<0$, and choosing $C_2$ large enough in terms of $c'$, the probability of having a set $S$ violating the statement is at most
\begin{align*}
    \sum_{k= C_2 (d\log d+\log n)}^\infty\exp\left(c'\frac{k}{d}\right)=\frac{\exp\left(d^{-1}c'C_2 (d\log d+\log n)\right)}{1-\exp\left(d^{-1}c'\right)}=o(1).
\end{align*}
\end{proof}

Lemma~\ref{key lemma} treats subsets of vertices from $V(L_1')$. The next step toward achieving the result of Theorem~\ref{setup_2_expansion} is to deal with subsets of vertices that have a non-empty intersection with $V(L_1)\setminus V(L_1')$ (Lemmas \ref{middle sets} and \ref{key lemma 2}). To that end, we prove an auxiliary lemma:

\begin{lemma}\label{close_L_1'}
    \textbf{Whp} $G$ contains at most $\frac{n}{\exp(\epsilon^{20} d)}$ vertices $v$ such that $|N_G(v)\cap V(L_1')|<\epsilon^{15} d$.
\end{lemma}

\begin{proof}
    We let $\delta=\epsilon^3$. Let $U_1$ be the set of vertices that participate in components of $G_2$ of size at least $\frac{1}{\delta}\log n$. Since our host graph $G$ satisfies the assumptions of Theorem~\ref{ERCP1}, due to Remark \ref{re: sprinkling_ercp1}, we observe that \textbf{whp} $U_1=V(L_1')$, so it is enough to show that \textbf{whp} $G$ contains at most $\frac{n}{\exp(\epsilon^{20} d)}$ vertices $v$ such that $|N_G(v)\cap U_1|<\epsilon^{15} d$.

    Let us define two more auxiliary sets:
    \begin{align*}
        U_2=\{v\in V(G)\mid |C_{G_2}(v)|\ge \delta d\}\quad \text{and} \quad U_3=\left\{v\in V(G)\mid \delta d\le |C_{G_2}(v)|\le \frac{1}{\delta}\log n\right\}.
    \end{align*}
    Note that $U_3\subseteq U_2$. Observe that we can apply Lemma~\ref{nointer_general} with $t\coloneqq \frac{1}{\delta}\log n$ to obtain that \textbf{whp} $|U_3|\le n\cdot \exp \left(-\epsilon^3\delta d\right)$. Since the graph $G$ is $d$-regular, and $d=\omega(1)$, we have trivially
    \[
        |N_G(U_3)|\le d\cdot n\cdot \exp \left(-\epsilon^3\delta d\right)\le n\cdot \exp \left(-\epsilon^4\delta d\right).
    \]
    
    Our goal is to bound from above the size of the set 
    \[
        B_1=\left\{w\in V(G)\mid |N_{G_2}(w)\cap U_1|<\epsilon^{15}d\right\}.
    \]
    Notice that by Lemma~\ref{W}, \textbf{whp} the set 
    \[
        B_2=\left\{w\in V(G)\mid |N_{G_2}(w)\cap U_2|<\epsilon^{15}d\right\}
    \]
    satisfies $|B_2|\le n\cdot \exp(-\delta^6 d)$. Moreover, observe that every vertex $v\in B_1$ is either also in $B_2$, or it has a neighbor in $U_3$, i.e. 
    \[
        B_1\subseteq B_2\cup N_G(U_3).
    \]
    Therefore, \textbf{whp}
    \[
        |B_1|\le |B_2|+|N_G(U_3)|\le \frac{n}{\exp\left(\epsilon^{20}d\right)},
    \]
    where we used $\delta=\epsilon^3$ and $d=\omega(1)$.
\end{proof}

We are ready to prove Lemma~\ref{middle sets} that deals with subsets of vertices outside $V(L_1')$.

\begin{lemma}\label{middle sets}
There exist constants $K\coloneqq K(\epsilon)>0$ and $\hat{c}\coloneqq \hat{c}(\epsilon)>0$ such that the following holds. Let $S\subseteq V(L_1')$ and $B\subseteq V(G)\setminus V(L_1')$ be such that $|S\cup B|\ge nd^{-C_1}$
and $|B|\ge K|S|$. Then, there exists a matching in $G_{p_3}$ of size at least $\hat{c}|B|$ between $B$ and $V(L_1')\setminus S$ with probability at least $1-\exp\left(-\hat{c}|B|\right)$.
\end{lemma}
\begin{proof}
By Lemma~\ref{close_L_1'}, \textbf{whp} $G$ contains at most $\frac{n}{\exp(\epsilon^{20} d)}$ vertices $v$ such that $|N_G(v)\cap V(L_1')|<\epsilon^{15} d$. In particular, \textbf{whp} every subset $M \subseteq V(G)$ contains at most $\frac{n}{\exp(\epsilon^{20} d)}$ vertices $v \in M$ such that $|N_G(v)\cap V(L_1')|<\epsilon^{15} d$. Under this assumption, we have 
\[
    e_G(B, L_1'\setminus S)\ge \left(|B|-\frac{n}{\exp(\epsilon^{20} d)}\right)\epsilon^{15} d-|S|d.
\]
Recall that we assume that $|S\cup B|\ge nd^{-C_1}$, that $|B|\ge K|S|$ and that $d$ is sufficiently large. Hence, $|B|\ge 2\cdot \frac{n}{\exp(\epsilon^{20} d)}$ and moreover there exists $c'\coloneqq c'(\epsilon,K)>0$ small enough such that
\[
    e_G(B, L_1'\setminus S)\ge \left(|B|-\frac{n}{\exp(\epsilon^{20} d)}\right)\epsilon^{15} d-|S|d\ge \frac{|B|}{2}\cdot \epsilon^{15} d-|S|d\ge c'|B|d.
\]
Thus, by Lemma~\ref{matching}, with probability at least $1-\exp\left(-c'\epsilon^6|B|\right)$ there exists a matching $M$ in $G_{p_3}$ of size at least $c'\epsilon^6|B|$ between $B$ and $L_1'\setminus S$.
\end{proof}

Equipped with Lemma~\ref{middle sets}, we show that sets in $L_1$, where the majority of their mass lies in $L_1\setminus L_1'$, admit good vertex expansion.

\begin{lemma}\label{key lemma 2}
There exist constants $K \coloneqq K(\epsilon)>0$ and $\hat{c}\coloneqq\hat{c}(\epsilon)>0$ such that \textbf{whp} for every $S_1\subseteq V(L_1')$ and $S_2\subseteq V(L_1)\setminus V(L_1')$ such that $G_p[S_1\cup S_2]$ is connected, $|S_2|\ge K|S_1|$, $|S_1\cup S_2|\ge nd^{-C_1}$, and $|S_1\cup S_2|\ge C_2 \log n$, we have
\begin{align*}
    |N_{L_1'}(S_1)|\ge \frac{\hat{c}|S_2|}{d},\quad \text{or}\quad 
    |N_{L_1}(S_1\cup S_2)|\ge \frac{\hat{c}|S_2|}{d}.
\end{align*}
\end{lemma}
\begin{proof}
We begin by exposing $G_2$. Now, given a set $S_1\subseteq V(L_1')$, we expose all the edges in $G_{p_3}$ that are either inside $V(G)\setminus V(L_1')$, inside $S_1$, or lie between $S_1$ and $V(G)\setminus V(L_1')$. Denote by $G_2'$ the graph $G_2$ together with these edges and by $L$ the giant component in $G_2'$ that contains $L_1'$. Note that $G_2 \subseteq G'_2 \subseteq G_3$ and $L_1'\subseteq L\subseteq L_1$. Further observe that for every $S_2\subseteq V(L_1)\setminus V(L_1')$ such that $S_1\cup S_2$ is connected in $L_1$, we have that $G_p[S_1\cup S_2]=G_2'[S_1\cup S_2]$.

Given $S_2\subseteq V(L_1)\setminus V(L_1')$ we denote $k_1\coloneqq |S_1|$ and $k_2\coloneqq |S_2|$. Let $K'\coloneqq K'(\epsilon)>0$ and $c'\coloneqq c'(\epsilon)>0$ be the constants given by Lemma~\ref{middle sets}, which will be needed later in the proof. We note that we may choose $K'$ sufficiently large. Set $K=K'$ and choose $\hat{c}\coloneqq \hat{c}(c')$ sufficiently small. Let us now facilitate a union bound argument. Let us further fix $k_2$ for now. We first count how many $S_1\subseteq V(L_1')$ with $|N_{L_1'}(S_1)|<\frac{\hat{c}k_2}{d}$. Given a set $S_1\subseteq V(L_1')$, suppose that $b_1\coloneqq |N_{L_1'}(S_1)|<\frac{\hat{c} k_2}{d}$, as otherwise the claim holds. Let us write $m$ for the number of components in $L_1'\setminus N_{L_1'}(S_1)$. Since $L_1'$ is connected and $G$ is $d$-regular, we have $m\le d\cdot b_1+1$. 

Hence, since $S_1$ has no neighbors in $L_1'$ outside $N_{L_1'}(S_1)$, it must be the union of some components of $L_1'\setminus N_{L_1'}(S_1)$, and so the number of ways to choose such an $S_1$ is at most $\binom{n}{b_1}2^{m}$. Thus, there are at most
\begin{align*}
    \sum_{b_1=1}^{\frac{\hat{c}k_2}{d}}\sum_{m=1}^{\hat{c}k_2}\binom{n}{b_1}2^{m}
    &\le \left(\frac{en}{\frac{\hat{c}k_2}{d}}\right)^{\frac{\hat{c}k_2}{d}}\cdot2^{\hat{c}k_2+1}\le \exp\left(\frac{\hat{c}k_2}{d}\left(\ln\left(\frac{end}{\hat{c}k_2}\right)+2d\right)\right)\le \exp\left(3\hat{c}k_2\right)
\end{align*}
sets $S_1\subseteq V(L_1')$ with $|N_{L_1'}(S_1)|<\frac{\hat{c}k_2}{d}$, where we used the assumptions that $|S_1\cup S_2|\ge nd^{-C_1}$ and that $|S_2|\ge K\cdot |S_1|$.

Now, given a set $S_1\subseteq V(L_1')$, we consider the number of ways to choose $S_2\subseteq V(L_1)\setminus V(L_1')$ that satisfy the conditions of the lemma, noting that having determined $S_1$, choosing $S_1\cup S_2$ determines $S_2$. We may assume that $b_2\coloneqq|N_{L}(S_1\cup S_2)|\le \frac{\hat{c}k_2}{d}$, since $L\subseteq L_1$. Since $L[S_1\cup S_2]$ is connected, and has all its neighbors in $N_L(S_1\cup S_2)$, exactly one of the at most $n$ components in $L\setminus N_L(S_1\cup S_2)$ is $S_1\cup S_2$. Since $S_1$ is fixed, we can identify this component. Hence, the number of ways to choose $S_1\cup S_2$ with $|N_L(S_1\cup S_2)|\le \frac{\hat{c}k_2}{d}$ is at most the number of ways to choose a set of size at most $\frac{\hat{c}k_2}{d}$ in $V(L)$. 
That is at most
\begin{align*}
   \sum_{b_2=1}^{\frac{\hat{c}k_2}{d}}\binom{n}{b_2}\le \left(\frac{en}{\frac{\hat{c}k_2}{d}}\right)^{\frac{\hat{c}k_2}{d}}\le \exp\left(\hat{c}k_2\right).
\end{align*}

We now turn to estimate the probability for the event $|N_{L_1}(S_1\cup S_2)|\le \frac{\hat{c}|S_2|}{d}$ to occur, for fixed $S_1\subseteq V(L_1')$ and $S_2\subseteq V(L_1)\setminus V(L_1')$ that satisfy the conditions of the lemma. 
Observe that for any $B\subseteq V(G)\setminus V(L_1')$, the event in Lemma~\ref{middle sets} depends only on edges in $G_{p_3}$ between 
$V(G)\setminus V(L_1')$ and $V(L_1')\setminus S_1$, which we have not yet exposed.

We now expose the edges between $V(G)\setminus V(L_1')$ and $V(L_1')\setminus S_1$, and use Lemma~\ref{middle sets} with $S_1$ in place of $S$ and $S_2$ in place of $B$ to obtain that the probability that $S_2$ has less than $c'|S_2|$ neighbors in $V(L_1')\setminus S_1$ in $G_{p_3}$ is at most $\exp\left(-c'|S_2|\right)$. Furthermore, the event that $S_2$ has at least $c'|S_2|$ neighbors in $V(L_1')\setminus S_1$ in $G_{p_3}$ clearly implies the lemma for $S=S_1\cup S_2$, since it implies that $|N_{L_1}(S_1\cup S_2)|\ge \frac{\hat{c}|S_2|}{d}$.

Therefore, for fixed $k_2$, the probability of an event violating the statement of the lemma is at most
\begin{align*}
    \exp\left(4\hat{c}k_2\right)\exp\left(-c'k_2\right)=\exp\left((4\hat{c}-c')k_2\right).
\end{align*}
Hence, since $\hat{c}$ is small enough in terms of $c'$, summing over all options of $k_2\ge \frac{C_2}{2}\log n$, we have, similarly to the proof of Lemma~\ref{key lemma}, that the probability of violating the statement of the lemma is $o(1)$ for sufficiently large $C_2$ in terms of $c'$.
\end{proof}

\begin{proof}[Proof of Theorem~\ref{setup_2_expansion}]
    Recalling that \textbf{whp} $L_1'\subseteq L_1$, we let $S_1=S\cap V(L_1')$ and $S_2=S\cap \left(V(L_1)\setminus V(L_1')\right)$. Let $c_{\ref{key lemma}}$ be the constant whose existence is asserted in Lemma~\ref{key lemma}, and let $K_{\ref{key lemma 2}}$ and $c_{\ref{key lemma 2}}$ be the constants whose existence is asserted in Lemma~\ref{key lemma 2}. Let $\beta>0$ be sufficiently small in terms of $c_{\ref{key lemma}}$, $K_{\ref{key lemma 2}}$ and $c_{\ref{key lemma 2}}$.
    Recall that we assume that $nd^{-C_1}\le |S|\le \frac{3\epsilon n}{2}$, that $|S|\ge C_2 (d\log d+\log n)$ and that $G_p[S]$ is connected.

   Suppose $|S_1|\ge \frac{|S_2|}{K_{\ref{key lemma 2}}}$. Then, $|S_1|\ge \frac{1}{K_{\ref{key lemma 2}}+1}\cdot nd^{-C_1}$ and $|S_1|\ge \frac{C_2}{K_{\ref{key lemma 2}}+1}(d\log d+\log n)$. By choosing $C_2$ sufficiently large in terms of $C_1,c_1,c_2,c_3,\alpha,\epsilon$, we can apply Lemma~\ref{key lemma} (with $t=\frac{1}{K_{\ref{key lemma 2}}+1}$) to obtain \textbf{whp} either 
   \begin{align*}
       |N_{L_1'}(S_1)|\ge \frac{c_{\ref{key lemma}}|S_1|}{d\log d}\ge\frac{\beta}{d\log d}|S|,
   \end{align*}
   or there is a family of at least $\frac{c_{\ref{key lemma}}|S_1|}{d}\ge \frac{\beta|S|}{d}$ many vertex-disjoint paths from $S_1$ to $V(L_1')\setminus S_1\subseteq V(L_1)\setminus S$. However, since each such path contributes a unique vertex to the neighborhood of $S$ in $V(L_1)$ (the first vertex along the path which is not in $S$), in the latter case $|N_{L_1}(S)|\ge \frac{\beta|S|}{d}$. As for the first case, $|N_{L_1}(S)|\ge |N_{L_1'}(S_1)|\ge\frac{\beta|S|}{d\log d}$, and so the result follows.

    Otherwise, $|S_2|\ge K_{\ref{key lemma 2}}|S_1|$ and so by Lemma~\ref{key lemma 2} \textbf{whp} either $|N_{L_1'}(S_1)|\ge \frac{c_{\ref{key lemma 2}}|S_2|}{d}$ or $|N_{L_1}(S)|\ge \frac{c_{\ref{key lemma 2}}|S_2|}{d}$. In the first case, $|N_{L_1}(S)|\ge |N_{L_1'}(S_1)|\ge \frac{\beta|S|}{d}$, and in the second case $|N_{L_1}(S)|\ge \frac{\beta|S|}{d}$. 
\end{proof}

\begin{proof}[Proof of Corollary~\ref{positive_diam}]
    We first show that \textbf{whp} for connected sets $S\subseteq V(L_1)$ we also have sufficiently good vertex expansion. Then, the typical good vertex expansion in $V(L_1)$ yield a control over the typical diameter of $L_1$.
    
    \begin{claim}\label{for_diam}
        There exists an absolute constant $c>0$ such that \textbf{whp} for any subset of the largest component $S\subseteq V(L_1)$ with $\frac{20\log n}{\epsilon^2}\le |S| \le nd^{-C_1}$ such that $S$ spans a connected subgraph in $G_p$, we have
    \begin{align*}
        \bigg|N_{G_p}(S)\bigg|\ge\frac{c\epsilon^2}{\log\left(\frac{1}{\epsilon}\right)}|S|.
    \end{align*}
    \end{claim}
    \begin{proof}[Proof of Claim~\ref{for_diam}]
        
    Consider the event $\mathcal{A}_k$,
    \begin{align*}
        \mathcal{A}_k=\left\{\exists S\subseteq V(G) \text{ such that } |S|=k, G_p[S]\text{ is connected} \ \& \  \big|N_{G_p}(S)\big|<\frac{\epsilon^2k}{40\log\left(\frac{1}{\epsilon}\right)}\right\}.
    \end{align*}
    We will show that for every $\frac{20\log n}{\epsilon^2}\le k\le nd^{-C_1}$, $\mathbb{P}\left(\mathcal{A}_k\right)=o\left(\frac{1}{n}\right)$, and therefore by the union bound over the $<n$ possible values of $k$, the probability of an event violating the statement of the theorem is $o(1)$.
    
    Let $S\subseteq V(G)$ with $|S|=k$ such that $G_p[S]$ is connected. Since it is connected, it must contain a spanning tree. By Claim~\ref{trees}, we have $n(ed)^{k-1}$ ways to choose a tree of order $k$, and we include its edges in $G_p$ with probability $p^{k-1}$. Now, consider the auxiliary random bipartite graph $\Gamma(S,p)$, whose one side is $S$, the other side is $N_G(S)$, and we include every edge of $G$ between $S$ and $N_G(S)$ in $\Gamma(S,p)$ with probability $p$. Then, we have that $|N_{G_p}(S)|\ge \nu\left(\Gamma(S,p)\right)$, where $\nu(H)$ is the matching number of a graph $H$. Thus, it suffices to bound the probability that a maximum matching in $\Gamma(s,p)$ is smaller than $\frac{\epsilon^2k}{40\log\left(\frac{1}{\epsilon}\right)}.$
    
    Let us first bound the probability that $\nu\left(\Gamma(S,p)\right)=i$. This is at most the probability that $\Gamma(s,p)$ has a maximal by inclusion matching of size $i$. We have at most $\binom{e_G\left(S,N_G(S)\right)}{i}\le \binom{kd}{i}$ ways to choose a matching $M$ of size $i$. We then need to include the edges of the matching, which happens with probability $p^i$. Due to the maximality of $M$, every edge of $G$ between $S$ and $N_G(S)$ disjoint from $M$ is not in $\Gamma(S,p)$. Thus, we have at least $e_G(S, S^c)-2id$ edges that do not fall into $\Gamma(S,p)$. Since $G$ satisfies $\mathrm{EE}\left(nd^{-C_1},\epsilon^3\right)$, we have
    \begin{align*}
       e_G\left(S, S^c\right)\ge \left(1-\epsilon^3\right)dk.
    \end{align*}
    Thus,
    \begin{align*}
        \mathbb{P}(\nu\left(\Gamma(S,p)\right)=i)\le \binom{kd}{i}p^i(1-p)^{\left(1-\epsilon^3\right)kd-2id}.
    \end{align*}
    Therefore, setting $u(\epsilon,k)\coloneqq \frac{\epsilon^2k}{40\log\left(\frac{1}{\epsilon}\right)}$, we obtain that
    \begin{align*}
        \mathbb{P}\left(\mathcal{A}_k\right)&\le n(ed)^{k-1}p^{k-1}\sum_{i=0}^{u(\epsilon,k)}\binom{kd}{i}p^i(1-p)^{\left(1-\epsilon^3\right)kd-2id}\\
        &\le n\left((1+\epsilon)\exp\left(1-(1+\epsilon)\left(1-\epsilon^3\right)\right)\right)^k\sum_{i=0}^{u(\epsilon,k)}\binom{kd}{i}p^i(1-p)^{-2id}.
    \end{align*}
    For every $i\ne 0$, we simplify the summands in the last summation as follows:
    \begin{align*}
        \binom{kd}{i}p^i(1-p)^{-2id}&\le \left(\frac{(1+\epsilon)ek}{i} \left(1-p\right)^{-2d}\right)^i\\
        &\le \left(\frac{(1+\epsilon)ek}{i} \exp \left(2d\left(\frac{1+\epsilon}{d}+\frac{(1+\epsilon)^2}{d^2}\right)\right)\right)^i\le \left(\frac{e^4k}{i}\right)^i,
    \end{align*}
    where we used the inequality $1-x\ge \exp\left(-x-x^2\right)$, which holds for any sufficiently small $x>0$. Thus,
    \begin{align*}
        \mathbb{P}\left(\mathcal{A}_k\right)\le n\left((1+\epsilon)\exp\left(1-(1+\epsilon)\left(1-\epsilon^3\right)\right)\right)^k\left(1+\sum_{i=1}^{u(\epsilon,k)}\left(\frac{e^4k}{i}\right)^i\right).
    \end{align*}
    Observe that the ratio of consecutive terms in the above sum is
    \begin{align*}
        \frac{\left(\frac{e^4k}{i}\right)^i}{\left(\frac{e^4k}{i+1}\right)^{i+1}}= \frac{(i+1)^{i+1}}{e^4ki^i}\le \frac{1}{2},
    \end{align*}
    and thus the sum is at most twice the final term. Hence,
    \begin{align*}
        \mathbb{P}\left(\mathcal{A}_k\right)&\le n\left((1+\epsilon)\exp\left(1-(1+\epsilon)\left(1-\epsilon^3\right)\right)\right)^k\left(1+2\left(\frac{e^440\log\left(\frac{1}{\epsilon}\right)}{\epsilon^2}\right)^{u(\epsilon,k)}\right)\\
        &\le 3n\left((1+\epsilon)\exp\left(-\epsilon+\epsilon^3+\frac{\epsilon^2}{40\log\left(\frac{1}{\epsilon}\right)}\log\left(\frac{41e^4\log\left(\frac{1}{\epsilon}\right)}{\epsilon^2}\right)\right)\right)^k\\
        &\le 3n\left((1+\epsilon)\exp\left(-\epsilon+\frac{\epsilon^2}{4}\right)\right)^k\le 3n\exp\left(-\frac{\epsilon^2k}{8}\right)=o\left(\frac{1}{n}\right),
    \end{align*}
    where the last inequality is due to $1+x\le \exp\left(x-\frac{3x^2}{8}\right)$ for $x$ small enough, and the equality is since we assumed $k\ge \frac{20\log n}{\epsilon^2}$. Union bound over the $<n$ different values of $k$ completes the proof.
    \end{proof}

Equipped with the expansion properties of connected subsets, we are now able to obtain a bound on the typical diameter of the largest component. By Theorem~\ref{L1}, \textbf{whp} we have that $|L_1|< 2\epsilon n$. Hence, if we can show that \textbf{whp} for any $v\in V(L_1)$, 
\begin{align}\label{eq: sufficient_for_diam}
    \Bigg|B_{L_1}\left(v,C(d\log^2 d+\log n)\right)\Bigg|\ge \epsilon n,
\end{align}
for some (large) constant $C>0$, then every two balls of such radius in $L_1$ intersect, and hence the typical diameter of the giant component of $G_p$ is $O(d\log^2 d+\log n)$.

Indeed, let $v\in V(L_1)$. By Claim~\ref{for_diam}, whenever $\frac{20\log n}{\epsilon^2}\le \big|B_{L_1}(v,r)\big| \le nd^{-C_1}$ for some $r\in\mathbb{N}$, we have \textbf{whp} $\big|N_{L_1}\left(B_{L_1}(v,r)\right)\big|\ge \frac{c\epsilon^2|B_{L_1}(v,r)|}{\log\left(\frac{1}{\epsilon}\right)}$. Also observe that $\big|B_{L_1}(v,r)\big|\ge \min\{r,|L_1|\}$ for any $r\in\mathbb{N}$. Therefore, for every $r\ge \frac{20\log n}{\epsilon^2}$ we have \textbf{whp}
\begin{align}\label{eq: ball_after_logn}
    \big|B_{L_1}(v, r+1)\big|\ge \min\left\{nd^{-C_1}, \left(1+\frac{c\epsilon^2}{\log\left(\frac{1}{\epsilon}\right)}\right)\big|B_{L_1}(v,r)\big|\right\}.
\end{align}
Consider
\begin{align*}
    r_1\coloneqq \frac{20\log n}{\epsilon^2}+\max\left\{0,\frac{\log n-C_1\log d}{\log \left(1+\frac{c\epsilon^2}{\log\left(\frac{1}{\epsilon}\right)}\right)}\right\}+C_2(d\log d+\log n)=O(d\log d+\log n).
\end{align*}
Observe that $|B_{L_1}(v,r_1)|\ge C_2 (d\log d+\log n)$, and by \eqref{eq: ball_after_logn}, $\left|B_{L_1}\left(v,r_1\right)\right|\ge nd^{-C_1}$. We may assume that $|B_{L_1}(v,r_1)|<\frac{3\epsilon n}{2}$, as otherwise we obtain \eqref{eq: sufficient_for_diam} as desired.

By Theorem~\ref{setup_2_expansion}, \textbf{whp} for any subset $S\subseteq V(L_1)$ with $nd^{-C_1}\le |S|\le \frac{3\epsilon n}{2}$ and $|S|\ge C_2 (d\log d+\log n)$, such that $G_p[S]$ is connected, we have that $\big|N_{G_p}(S)\big|\ge \frac{\beta}{d\log d}|S|$. Therefore, for any $r\ge r_1$, we have that \textbf{whp}
\begin{align*}
    \big|B_{L_1}(v,r+1)\big|\ge\min\left\{\frac{3\epsilon n}{2}, \left(1+\frac{\beta}{d\log d}\right)\big|B_{L_1}(v,r)\big|\right\}.
\end{align*}
Since $\big|B_{L_1}(v,r_1)\big|\ge nd^{-C_1}$, for $r_2\coloneqq \frac{\log\left(\frac{3\epsilon}{2}d^{C_1}\right)}{\log\left(1+\frac{\beta}{d\log d}\right)}$ we have that \textbf{whp}
\begin{align*}
    \Bigg|B_{L_1}\left(v,r_1+r_2\right)\Bigg|\ge \frac{3\epsilon n}{2}> \frac{|V(L_1)|}{2}.
\end{align*}
Thus, \textbf{whp} $\mathrm{diam}(L_1)\le 2(r_1+r_2)=O(d\log^2 d+\log n)$.
\end{proof}

\section{Diameter of the Giant: Upper Bound --- Proof of Theorem~\ref{pure_diam}}\label{Ppure_diam}

Throughout this section we use notations and assumptions in Theorem~\ref{pure_diam}. We utilize a double-exposure/sprinkling argument. Define $p_2=\frac{\epsilon^3}{d}$ and let $p_1$ be such that $(1-p_1)(1-p_2)=1-p$, noting that $p_1\geq \frac{1+\epsilon-\epsilon^3}{d}$. Notice that $G_p$ has the same distribution as $G_{p_1}\cup G_{p_2}$. We abbreviate $G_1\coloneqq G_{p_1}$ and $G_2\coloneqq G_{p_1}\cup G_{p_2}$. Denote by $L_1'$ and $L_1$ the largest components of $G_{p_1}$ and $G_p$, respectively. Note that \textbf{whp} the size of the second largest component in $G_2$ is strictly smaller than $|V(L_1')|$, and hence \textbf{whp} $L_1'\subseteq L_1$. 

Let us first discuss the proof’s strategy: we show that \textbf{whp} for every $v\in V(L_1')$ we have $|B_{G_2}(v,r)\cap V(L_1')|>\frac{|V(L_1')|}{2}$, for some $r=O\left(d\log n\left(1+\log_d\log n\right)\right)$. Thus, \textbf{whp} the distance (in $G_2$) between any two vertices in $L_1'$ is $O\left(d\log n\left(1+\log_d\log n\right)\right)$. We also show that \textbf{whp} for every $v\in V(L_1)$, the distance (in $G_2$) between $v$ and $L_1'$ is $O(\log n)$. Combining these claims, we obtain \textbf{whp} $\mathrm{diam}(L_1)=O\left(d\log n\left(1+\log_d\log n\right)\right)$.

We start by showing that \textbf{whp} the distance (in $G_2$) between every $v\in L_1\setminus L_1'$ and $L_1'$ is $O(\log n)$. This is established by the following argument, which generalizes Lemma 3.2 in \cite{zbMATH07628800}.

\begin{lemma}\label{out_L_1'}
    Given a vertex $v\in V(L_1')$, denote by $C_v$ the set of vertices contained in components of $L_1\setminus L_1'$ in $G_2$ that are adjacent to $v$ in $G_2$. Then there exists a constant $K \coloneqq K(\epsilon)>0$ such that \textbf{whp} $|C_v| \le K\log n$ for every $v\in V(L_1')$.
\end{lemma}

\begin{proof}
We first note that by Theorem~\ref{ERCP1}, there exists a constant $K'\coloneqq K'(\epsilon)$ such that \textbf{whp} any component of $G_1$ except $L'_1$ has order at most $K'\log n$. Let $K \ge \frac{10}{\epsilon^2}$.

Suppose that there is some vertex $v \in V(L'_1)$ such that $|C_v|\geq K \log n$. We note that $C_v \cup \{v\}$ is connected in $G_2$, and $C_v$ is the disjoint union of some sets $\{C_1,\ldots, C_r\}$ where each $C_i$ is the vertex set of some component of $G_1$, each of which has size at most $K' \log n$. It follows that there must be some subset $\hat{C} \subseteq C_v$ such that $\hat{C} \cup \{v\}$ is connected in $G_2$, $K \log n \leq |\hat{C}| \leq (K' + K)\log n$ and $\hat{C}$ is the union of some subset of $\{C_1,C_2, \ldots, C_r\}$.

In particular, there is some spanning tree $T$ of $\hat{C} \cup \{v\}$, all of whose edges are present in $G_2$, such that no edge in the edge-boundary of $V(T) \setminus \{v\}$ is present in $G_1$. 

Let us bound the probability that such a tree of order $k$ exists in $G_2$ for each 
\[
    K\log n +1 \leq k \leq (K' + K)\log n+1.
\] 
Set $I\coloneqq [ K \log n +1,(K' + K)\log n+1]$. We can fix such a tree $T$ by choosing a root vertex $v$ and choosing one of the at most $(ed)^{k-1}$ possible rooted trees of order $k$ with root $v$ in $G$, where we have bounded the number of possible trees by Claim~\ref{trees}.

Now, $T$ has $k-1$ edges, and by \ref{prop2P}, $e_G(T\setminus\{v\},(T\setminus \{v\})^c)\ge \left(1-\epsilon^3\right)d(k-1)$. Note that each edge is in $G_2$ with probability $p$ and each edge is not in $G_1$ with probability $1-p_1$, and that, while these two events are not necessarily independent, they are clearly negatively correlated. 

It follows by the union bound that the probability that such a tree of order $k$ exists in $G_2$ is at most

\begin{equation*}
n(ed)^{k-1} p^{k-1} \left(1- p_1\right)^{\left(1-\epsilon^3\right)d(k-1)}.
\end{equation*}
In particular, the probability that such a tree exists for $k \in I$ is at most
\begin{align*}
&n\sum_{k\in I} (e(1+\epsilon))^{k-1}\left(1- \frac{1+\epsilon-\epsilon^3}{d}\right)^{\left(1-\epsilon^3\right)d(k-1)}\\
&\leq  n\sum_{k\in I} \exp\left( (k-1)\left( 1 + \log(1+\epsilon) - (1-\epsilon^3)(1+\epsilon-\epsilon^3) \right)\right)\\
&\leq n \sum_{k\in I} \exp\left( -\frac{(k-1)\epsilon^2}{5}\right) = o(1),
\end{align*}
where the last inequality holds since $\log (1+ x) \le x-\frac{x^2}{3}$ for small enough $x>0$ and since $K \ge \frac{10}{\epsilon^2}$.
\end{proof}

We now turn to bounding the typical distance (in $G_2$) between vertices of $L_1'$.

\begin{lemma}\label{diam_L_1'}
    There is an integer $r=O\left(d\log n\left(1+\log_d\log n\right)\right)$ such that \textbf{whp} for every $v\in V(L_1')$, $|B_{G_2}(v,r)\cap V(L_1')|>\frac{|V(L_1')|}{2}$.
\end{lemma}

\begin{proof}
    We begin by exposing $G_1$. Let $v\in V(L_1')$, and let $c_{\ref{paths_setup2}}$ be the constant whose existence is asserted in Lemma~\ref{paths_setup2}. Now, expose the edges of $G_{p_2}$. Define $A_0=B_{G_1}\left(v,\frac{3}{c_{\ref{paths_setup2}}}d\log n\right)$ and $B_0=V(L_1')\setminus A_0$. We may assume that $|A_0|\le |B_0|$, as otherwise we are done. By Lemma~\ref{paths_setup2}, with probability $1-\exp\left(-\frac{c_{\ref{paths_setup2}}|A_0|}{d}\right)$, there exists a family of $\frac{c_{\ref{paths_setup2}}|A_0|}{d}$ vertex-disjoint $A_0-B_0$-paths of length at most $2(\log_d\log n+1)+1$ in $G_{p_2}$. Define $A_1$ to be $A_0$ with the endpoints of the at least $\frac{c_{\ref{paths_setup2}}|A_0|}{d}$ paths between $A_0$ and $B_0$. Also, set $B_1=V(L_1')\setminus A_1$. 
    Observe the following facts: first, $|A_1|\ge |A_0|+\frac{c_{\ref{paths_setup2}}|A_0|}{d}=|A_0|\left(1+\frac{c_{\ref{paths_setup2}}}{d}\right)$; second, $A_1\subseteq B_{G_2}\left(v,\frac{3}{c_{\ref{paths_setup2}}}d\log n+2(\log_d\log n+1)+1\right)\cap V(L_1')$. 
    
    As long as $|A_k|\le \frac{|V(L_1')|}{2}$, with probability $1-\sum_{i=0}^{k+1}\exp\left(-\frac{c_{\ref{paths_setup2}}|A_i|}{d}\right)$, we can define $A_{k+1}$ as $A_k$ together with the endpoints of the at least $\frac{c_{\ref{paths_setup2}}|A_k|}{d}$ paths between $A_k$ and $B_k$ and $|A_{k+1}|\ge |A_k|\left(1+\frac{c_{\ref{paths_setup2}}}{d}\right)$ and that $A_{k+1}\subseteq B_{G_2}\left(v,\frac{3}{c_{\ref{paths_setup2}}}d\log n+(k+1)(2(\log_d\log n+1)+1)\right)\cap V(L_1')$. We also set $B_{k+1}=V(L_1')\setminus A_{k+1}$. Denote by $r_0=\log_{1+\frac{c_{\ref{paths_setup2}}}{d}}\left(\frac{|V(L_1')|}{2}\right)=O(d\log n)$. Notice that with probability at least
    \begin{align*}
        1-\sum_{i=0}^{r_0+1}\exp\left(-\frac{c_{\ref{paths_setup2}}|A_i|}{d}\right)\ge 1-n\cdot \exp\left(-\frac{c_{\ref{paths_setup2}}|A_0|}{d}\right)=1-o\left(\frac{1}{n}\right)
    \end{align*}
    we have $A_{r_0+1}$ such that
    \begin{align*}
        |A_{r_0+1}|\ge |A_0|\left(1+\frac{c_{\ref{paths_setup2}}}{d}\right)^{r_0+1}>\frac{|V(L_1')|}{2}
    \end{align*}
    and furthermore
    \begin{align*}
        A_{r_0+1}\subseteq B_{G_2}\left(v,\frac{3}{c_{\ref{paths_setup2}}}d\log n+(r_0+1)\cdot (2(\log_d\log n+1)+1)\right)\cap V(L_1')\subseteq B_{G_2}(v,r)\cap V(L_1')
    \end{align*}
    for $r=\frac{3}{c_{\ref{paths_setup2}}}d\log n+(r_0+1)\cdot (2(\log_d\log n+1)+1)=O\left(d\log n\left(1+\log_d\log n\right)\right)$. To end, by the union bound over all vertices $v\in V(L_1')$ we have that \textbf{whp} $|B_{G_2}(v,r)\cap V(L_1')|>\frac{|V(L_1')|}{2}$ for every $v\in V(L_1')$.
\end{proof}

Lemmas \ref{out_L_1'} and \ref{diam_L_1'} yield that \textbf{whp} the distance between any pair of vertices of $L_1$ is at most $O(\log n)+O\left(d\log n\left(1+\log_d\log n\right)\right)=O\left(d\log n\left(1+\log_d\log n\right)\right)$.

\section{No Giant: Tightness of Theorem~\ref{L1} --- Proof of Theorem~\ref{no_L1}}\label{pno_L1}

Let $\sigma>1,\alpha>0$ be some constants. Set $C>1$ such that $b=\frac{\alpha}{60\sigma e^2C}$. Let $d, d_1$, and $n$ be even positive integers satisfying the following conditions:
\begin{align*}
    d=\omega(1);\quad d_1= \frac{d}{\sigma e^2}; \quad n=\omega(d^2).
\end{align*}

Finally, set $p=\frac{\sigma}{d}$. We construct a $d$-regular $n$-vertex graph $G$ as follows.

We start with $n_1\coloneqq \frac{n}{\frac{4}{\alpha}Cd+1}$ copies of a $\left(\frac{4}{\alpha}Cd,d,\mu\right)$-graph, with $\mu=O\left(\sqrt{d}\right)$ (notice that such a graph exists since $\frac{4}{\alpha}Cd\cdot d$ is even, and thus a random $d$-regular graph on $\frac{4}{\alpha}Cd$ vertices satisfies this --- see, for example, \cite[Theorem 1.1]{zbMATH07751065}. See also \cite{zbMATH05037081,zbMATH07036340} and the references therein for even stronger results).
For the $i$-th copy, $i=1,\ldots,n_1$, we find a matching $E_i$ of size $\frac{d-d_1}{2}$ in the graph (such a matching exists by \cite[Theorem 4.3]{zbMATH05037081}, for example), and remove it. We denote by $F_i$ the set of the $d-d_1$ vertices participating in the matching $E_i$. Finally, we set $H_i$ as the $i$-th copy minus the matching $E_i$.

Now, let $H_0$ be an $\left(n_1,d_1,\lambda\right)$-graph, with $\lambda=O\left(\sqrt{d_1}\right)$ (notice that such a graph exists by the same argument as before and since $n_1\cdot d_1=\frac{n}{\frac{4}{\alpha}Cd+1}\cdot d_1$ is even). Label the vertices of $H_0$ by $V(H_0)=\{v_i\}_{i=1}^{n_1}$. 

Our desired graph $G$ will be the graph $H_0$, such that for every $1\le i\le n_1$, we connect $v_i$ with $H_i$ via all vertices in $F_i$.

Observe that $G$ is a $d$-regular graph on $n_1\cdot \left(\frac{4}{\alpha}Cd+1\right)=n$ vertices satisfying properties $\mathrm{iso}(b)$ and $\mathrm{EE}(Cd,\alpha)$. It is not hard to verify that $G$ satisfies these two properties, and its proof is placed in Appendix \ref{app: 5} since the verification is somewhat tedious.

All that remains to prove is that \textbf{whp} all components of $G_p$ are of size at most $\frac{9}{\alpha}Cd\log n$. 
Indeed, let $k=\frac{9}{\alpha}Cd\log n$, and assume $G_p$ has a component of size at least $k$. In particular, there is a tree $T$ of order $k$ such that $T\subseteq G_p$. Observe that, by the construction of $G$, there must be a subtree $T'\subseteq T$ of order $k_0\coloneqq \frac{k}{\frac{4}{\alpha}Cd+1}\ge 2\log n$ such that $T'\subseteq (H_0)_p$. We bound from above the probability that such an event occurs. By Claim~\ref{trees}, there are at most $n_1(ed_1)^{k_0-1}$ ways to choose a tree on $k_0$ vertices in $H_0$, and we retain its edges in $(H_0)_p$ with
probability $p^{k_0-1}$. Therefore, the probability that $(H_0)_p$ has a tree of order at least $k_0$ is at most
\begin{align*}
    \frac{n}{\frac{4}{\alpha}Cd+1}(ed_1)^{k_0-1}p^{k_0-1}=\frac{n}{\frac{4}{\alpha}Cd+1}\left(\frac{1}{e}\right)^{\frac{k}{\frac{4}{\alpha}Cd+1}-1}=o(1),
\end{align*}
where we used $pd_1=\frac{1}{e^2}$.

\section{Lack of Vertex Expansion --- Proof of Theorem~\ref{no_expansion}}\label{Pno_expansion}

Let $\alpha>0$ and $0<c_1<\frac{1}{160\epsilon}$ and set $C=\frac{1}{10\epsilon}$. Assume that both $d$ and $C$ are even. We construct a $d$-regular $n$-vertex graph $G$ as follows: for every $1\le i\le \frac{n}{2d\log n}\eqqcolon k$, we consider two copies of $\{1,2,\ldots,d\log n\}$ and denote them by $A_i$ and $B_i$. We let $\Gamma_1$ be a random $(d-C)$-regular graph on $d\log n$ vertices and $\Gamma_2$ be a random $C$-regular bipartite graph on $d\log n$ vertices each side.
\begin{enumerate} [(\arabic*)]
    \item For each $1\le i\le k$, we set $H^i\coloneqq G[A_i]$ as a copy of $\Gamma_1$.
    \item For each $1\le i\le k$, we set $G[A_i,B_i]$ as a copy of $\Gamma_2$.
    \item Finally, we set $H^0\coloneqq G[\cup_{i=1}^k B_i]$ as a random $(d-C)$-regular graph.
\end{enumerate} 

Observe that since $C$ is a sufficiently large constant, for every $1\le i\le k$, due to \cite[Theorem 4]{zbMATH07671057}, \textbf{whp} we have $\lambda\coloneqq\lambda_2\left(G[A_i,B_i]\right)=O\left(\sqrt{C}\right)$ and we further notice that \textbf{whp} $H^i$ is a $\left(d\log n,d-C,\nu\right)$-graph, where $\nu=O\left(\sqrt{d-C}\right)$. Finally, notice that \textbf{whp} $H^0$ is a $\left(\frac{n}{2},d-C,\tau\right)$-graph for $\tau=O\left(\sqrt{d-C}\right)$ (see, e.g., \cite[Theorem 1.1]{zbMATH07751065}. See also \cite{zbMATH05037081,zbMATH07036340} and the references therein for even stronger results). By construction, $G$ is a $d$-regular $n$-vertex graph. Proofs that the graph satisfies $\mathrm{iso}(c_1)$, $\mathrm{VE}(c_2\log n,\alpha)$ and $\mathrm{EE}(c_2 d\log n,\alpha)$ can be found in Appendix \ref{app: 6}.

Set $p=\frac{1+\epsilon}{d}$. We turn to proving that \textbf{whp} the largest component $L_1$ in $G_p$ contains a subset $S$ such that $|S|=\Theta(d\log n)$, $G_p[S]$ is connected, and
\begin{align*}
    |N_{G_p}(S)|\le \frac{\beta}{d\log n}|S|,
\end{align*}
for some $\beta\coloneqq \beta(\epsilon)>0$. In other words, we show that \textbf{whp} there is a subset $S\subseteq V(L_1)$ as above, such that $|N_{G_p}(S)|=\Theta_\epsilon(1)$.

In this proof, we expose the edges of $G$ and $G_p$ in the following order. At first, we expose the edges of $H^i$ for all $0\le i\le k$. Then apply $p$-bond percolation on $H^0$ and afterwards on $H^i$, for all $1\le i\le k$. Now, expose the edges of $G[A_i,B_i]$ for all $1\le i\le k$, and finally perform $p$-bond percolation on them. The order in which the edges are exposed does not matter, since these sets of edges are independent of each other.

Since $d=\omega(1)$, it holds that $\sigma_0\coloneqq p\cdot (d-C)>1$. Apply Lemma~\ref{FKM} with $\sigma_0$ to obtain that \textbf{whp} there exists a unique giant component $K\coloneqq L_1(H^0)$ in $H^0_p$ of size at least $\epsilon|V(H^0)|=\frac{\epsilon n}{2}$. Notice that our graph $G$ satisfies the assumptions of Theorem~\ref{L1}, hence $K$ is part of the unique giant component of $G_p$. Notice that there are at least $\frac{\epsilon n}{4d\log n}$ many indices $1\le i\le k$ for which $|V(K)\cap B_i|\ge \frac{\epsilon}{4}\cdot d\log n$, as otherwise the number of vertices in $K$ is at most 
\begin{align*}
    \left(\frac{\epsilon n}{4d\log n}-1\right)\cdot |B_i|+\left(k-\frac{\epsilon n}{4d\log n}+1\right)\cdot \frac{\epsilon}{4}\cdot d\log n\le \frac{\epsilon n}{4}+\frac{\epsilon n}{8}<\frac{\epsilon n}{2},
\end{align*}
which yields a contradiction. Denote by $I'\subseteq [k]$ the set of indices $i$ such that $|V(K)\cap B_i|\ge \frac{\epsilon}{4}\cdot d\log n$, and set $V_i\coloneqq V(K)\cap B_i$. The above argument showed that $|I'|\ge \frac{\epsilon n}{4d\log n}$. For each $i\in I'$, consider $H^i$. We denote by $U_i$ the vertex set of the largest connected component of $H^i_p$. Since $\sigma_1\coloneqq p\cdot (d-C)>1$, by Lemma~\ref{FKM} with $\sigma_1$, we obtain that \textbf{whp}
\[
    \frac{\epsilon}{2}\cdot d\log n\le |U_i|\le 2\epsilon d\log n.
\]

Denote by $I\subseteq I'$ the set of indices $i$ such that $\frac{\epsilon}{2}\cdot d\log n\le |U_i|\le 2\epsilon d\log n$. Since the subgraphs $\{H^i\}_{i\in I'}$ are disjoint, using a standard Chernoff bound implies that \textbf{whp} $|I|\ge \frac{|I'|}{2}\ge \frac{\epsilon n}{8d\log n}$. We continue assuming that this holds deterministically. Summing up, we have that there are at least $\frac{\epsilon n}{8d\log n}$ many indices $i$ for which
\begin{align}\label{size_U,V}
    \frac{\epsilon}{2}\cdot d\log n\le |U_i|\le 2\epsilon d\log n \quad \text{and} \quad |V_i|\ge \frac{\epsilon}{4}\cdot d\log n. 
\end{align}

Since $V_i\subseteq K$, by proving that \textbf{whp} there is at least one $i\in I$ for which $e_{G_p}(U_i,V_i)=1$ and that $e_{G_p}(U_i,B_i\setminus V_i)=0$, we yield that $U_i$ is part of a giant component of $G_p$, such that $|U_i|\ge \frac{\epsilon}{2}d\log n$ and also $|N_{G_p}(U_i)|=1\le \frac{2}{\epsilon\cdot d\log n}\cdot |U_i|$, which establishes the theorem.

Indeed, let $i\in I$. Let $X_i$ be the indicator of the event that $e_{G_p}(U_i,V_i)=1$ and that $e_{G_p}(U_i,B_i\setminus V_i)=0$. In other words, $X_i$ is the indicator of the event that from all $C|U_i|$ edges between $U_i$ and $B_i$, exactly one of them falls into $G_p$, and this edge is an edge between $U_i$ and $V_i$. Since $G[A_i,B_i]$ is a random $C$-regular graph, then using the configuration model, we can couple the distribution of $e_{G[A_i,B_i]}(U_i,V_i)$ with $\mathrm{Hypergeom}(Cd\log n, C|V_i|, C|U_i|)$. Hence, by Claim~\ref{cher_hyp}, \textbf{whp}
\[
    e_{G[A_i,B_i]}(U_i,V_i)\ge \frac{C|U_i||V_i|}{2d\log n}.
\]
Therefore, using \eqref{size_U,V}, we obtain
\begin{align*}
        \mathbb{P}(X_i=1)&\ge \frac{C|U_i||V_i|}{2d\log n}\cdot p\left(1-p\right)^{C|U_i|-1}\ge \frac{C}{2}\cdot \frac{\epsilon}{2}\cdot d\log n\cdot \frac{\epsilon}{4}\cdot \frac{1+\epsilon}{d}\cdot \left(1-\frac{1+\epsilon}{d}\right)^{C\cdot 2\epsilon d\log n}\\
        &\ge \exp\left(-\left(\frac{1+\epsilon}{d}+\frac{(1+\epsilon)^2}{d^2}\right)\cdot C\cdot 2\epsilon d\log n\right)\ge n^{-4\epsilon C},
    \end{align*}
    where the penultimate inequality holds since $1-x\ge \exp\left(-x-x^2\right)$ for sufficiently small $x>0$.
    Let $X=\sum_{i\in I}X_i$. By the above calculations, we have
    \begin{align*}
        \mathbb{E}[X]\ge |I|\cdot n^{-4\epsilon C}\ge \frac{\epsilon n}{8d\log n}\cdot n^{-4\epsilon C}=\frac{\epsilon n^{1-4\epsilon C}}{8 d\log n}.
    \end{align*}
    Hence, by Claim~\ref{cher}, \textbf{whp} $X\ge \frac{\epsilon n^{1-4\epsilon C}}{8 d\log n}$. We assumed that $C=\frac{1}{10\epsilon}$ and that $d<\sqrt{n}$ and hence \textbf{whp} $X\ge 1$.

\section{Diameter of the Giant: Lower Bound --- Proof of Theorem~\ref{general_large_diam}}\label{P_larege_diam}

Let $c_1,\alpha>0$ and $\sigma>1$ be constants and choose $c_2\coloneqq c_2(\alpha)>0$ to be small enough constant. Let $n,d$ be sufficiently large integers such that $d=o\left(n\right), d=\omega(1)$. Assume further that $d$ is divisible by $4$. Set $p\ge\frac{\sigma}{d}$. Let $t\ge \frac{4c_1}{\alpha}\cdot\frac{n}{d}$ and also $t\ge 2d$. Let us further assume that $t<\frac{n}{2}$, as otherwise the conclusion of Theorem~\ref{general_large_diam} is trivial.

Let $\Gamma_1$ be a random $\frac{\alpha d}{4}$-regular bipartite graph with sides of size $t$ each. Due to \cite[Theorem 4]{zbMATH07671057}, \textbf{whp} $\lambda_2(\Gamma_1)= O\left(\sqrt{d}\right)$. Also, let $\Gamma_2$ be a random $\left(1-\frac{\alpha}{2}\right)d$-regular graph on $t$ vertices. Notice that $\Gamma_2$ is a $\left(t,\left(1-\frac{\alpha}{2}\right)d,\tau\right)$-graph for $\tau=O\left(\sqrt{d}\right)$ (see, e.g., \cite[Theorem 1.1]{zbMATH07751065}. See also \cite{zbMATH05037081,zbMATH07036340} and the references therein for even stronger results).

Consider $\frac{n}{t}$ sets of $t$ vertices each. We denote the $i$-th set by $S_i$. For the rest of this section, all indices related to $\{S_i\}_{i=1}^{n/t}$ are taken modulo $n/t$. In particular, we identify $S_{\frac{n}{t}+1}$ with $S_1$.

Our desired graph $G$ is constructed as follows: $V(G)=\bigcup_{i=1}^{n/t} S_i$ such that for every $1\le i\le \frac{n}{t}$, embed a copy of $\Gamma_1$ in $G[S_i,S_{i+1}]$ and embed a copy of $\Gamma_2$ in $G[S_i]$. Note that since the indices are taken modulo $n/t$, there is a bipartite graph between $S_{\frac{n}{t}}$ and $S_1$. For every $1\le i\le \frac{n}{t}$, denote $H_i \coloneqq G[S_i,S_{i+1}]$ and $F_i\coloneqq G[S_i]$. We note that $G$ is a $d$-regular $n$-vertex graph.

Now, observe that every connected subgraph $S$ of $G_p$ must be spread over at least $|S|/t$ slices of $\{S_i\}_{i=1}^{n/t}$. In particular, the giant component $L_1\subseteq G_p$ has diameter at least $\frac{|L_1|}{2t}-1$.

It remains to show that $G$ satisfies $\mathrm{iso}(c_1)$, $\mathrm{VE}(c_2\cdot \frac{t}{d},\alpha)$ and $\mathrm{EE}(\frac{\alpha}{2}t,\alpha)$. The proof are somewhat tedious, but straightforward, and are given in Appendix \ref{app: 7}. Let us note that the first property $\mathrm{iso}(c_1)$ holds due to the edges between the slices $\{S_i\}_{i=1}^{n/t}$, whereas the other two properties $\mathrm{VE}(c_2\cdot \frac{t}{d},\alpha)$ and $\mathrm{EE}(\frac{\alpha}{2}t,\alpha)$ follow from the edges inside the slices. We also observe that the trade-off between the sizes of the sets and the quality of the expansion in property $\mathrm{EE}(\frac{\alpha}{2}t,\alpha)$ stems from the fact that achieving almost optimal edge expansion (a factor of nearly $d$) requires controlling the sizes of the sets within each graph $F_i=G[S_i]$ ($1\le i\le \frac{n}{t}$). For further details, see Appendix \ref{app: 7}.

\section{Discussion}\label{disc}

In this paper, we developed a unified approach to bond percolated graphs obtained from host graphs with three different sets of vertex and edge expansion properties. This approach both yields a variety of new results and brings together some existing ideas and techniques, highlighting their generality and applicability. Our results include a new (tight) scenario for which a unique giant component of the percolated graph is typically guaranteed and the other components are relatively small. We further provided results concerning typical vertex expansion of the giant component and several combinatorial properties and aspects of it, in particular the diameter of the giant component of the percolated graph, for which we established both upper and lower bounds.

Theorems \ref{general} and \ref{setup_2_expansion} show a typical vertex expansion of the giant component obtained by supercritical bond percolation on host graphs with fairly good vertex and edge expansion properties. These results enable us to find \textbf{whp} a cycle of length $\Omega\left(\frac{n}{d\log d}\right)$ in the giant component, as well as to obtain a lower bound of order $\frac{\sqrt{n}}{d\log d}$ on the typical Hadwiger number of the giant component. Both results are optimal up to polynomial terms in $d$. 

We further note that the Hadwiger number $h(H)$ of a graph $H$ is related to the genus of the graph --- a topological property of the graph that generalizes planarity. Formally, the genus of a graph $H$ is the minimum integer $g\coloneqq g(H)$ such that $H$ has an embedding into an orientable surface of genus $g$ such that there are no crossing edges. The relation between $h(H)$ and $g(H)$ is given by
\[
    g(H)\ge \left\lfloor\frac{(h(H)-3)(h(H)-4)}{12}\right\rfloor.
\]
Therefore, our result on the Hadwiger number of the giant component (Corollary~\ref{properties_general} \ref{cor_2}) implies that the genus of $L_1$ is at least linear in the number of vertices, up to polynomial terms in $d$. 

Theorem~\ref{setup_2_expansion} has another implication regarding the diameter of the giant component: Corollary~\ref{positive_diam} shows that assuming $\mathrm{iso}(c_1)$, $\mathrm{VE}(c_2\log n,1-c_3)$ and $\mathrm{EE}(\max\left\{Cd\log n,nd^{-C_1}\right\},\epsilon^3)$ on the host graph $G$, for some appropriate parameters $c_1,c_2,c_3,C,C_1$, we obtain that the giant component $L_1$ satisfies \textbf{whp} 
\[
    \mathrm{diam}(L_1)=O\left(d\log^2 d+\log n\right).
\]

We have also established an upper bound on the typical diameter of the giant component where the host graph admits only the (by now) standard assumptions that guarantee a unique giant component of the percolated graph (Theorem~\ref{ERCP1}). More formally, assuming $\mathrm{iso}(c_1)$, $\mathrm{VE}(c_2\log n,1-c_3)$ and $\mathrm{EE}(Cd\log n,\epsilon^3)$ on the host graph $G$, for some appropriate parameters $c_1,c_2,c_3,C$, we obtain that the giant component $L_1$ satisfies \textbf{whp} 
\[
    \mathrm{diam}(L_1)=O\left(d\log n(1+\log_d\log n)\right).
\]

We now turn to the results in the other direction, concerning lower bounds on the diameter. We have found a significant gap between well-studied graphs, for example, random subgraphs of $K_n$, $(n,d,\lambda)$-graphs, and graphs that satisfy fairly good vertex and edge expansion properties, as can be seen in Theorems \ref{isop}, \ref{ERCP1}, or \ref{L1}. In particular, Corollaries \ref{large_diam_setup1} and \ref{large_diam_setup2} show that although the assumptions in Theorems \ref{isop}, \ref{ERCP1}, or \ref{L1} are sufficient to guarantee \textbf{whp} the emergence of the giant component of the percolated graph, there are graphs that satisfy the assumptions of Theorems \ref{isop}, \ref{ERCP1}, or \ref{L1}, but the giant component $L_1$ satisfies \textbf{whp}
\[
    \mathrm{diam}(L_1)=\Omega\left(\min\left\{d,\frac{n}{d}\right\}\right)\quad \text{or}\quad \mathrm{diam}(L_1)=\Omega\left(\min\left\{d,\frac{n}{d\log n}\right\}\right).
\]
Therefore, the conditions of Theorems \ref{isop}, \ref{ERCP1}, or \ref{L1} are sufficient to guarantee \textbf{whp} the emergence of the giant component of the percolated graph (where Theorem~\ref{ERCP1} also ensures the optimal upper bound on the second-largest component of the percolated graph), yet they are not sufficient to establish that \textbf{whp} the diameter of the giant component is $O(\log n)$ which hold when the host graphs are $K_n$ or $(n,d,\lambda)$-graphs. It is worth noting that we do obtain a fairly tight result under the assumptions of Theorem~\ref{ERCP1} concerning the diameter of the giant component. Indeed, assuming that $d<\sqrt{\frac{n}{\log n}}$, the lower and upper bounds differ only by a factor of $\log n(1+\log_d\log n)$. Therefore, it would be very interesting to explore what the exact asymptotic behavior of the diameter of the giant component is, if it exists, under this regime.

This is a first natural direction for future research --- improving the results presented in this paper and obtaining tighter bounds, as well as to determining under what assumptions on the host graph, the results presented in this paper are tight. For example, as mentioned above, we have established upper and lower bounds on the typical diameter of the giant component whose dominant terms are $d$ and $\log n$. Taking the recent breakthrough result on the diameter of the giant component of the percolated hypercube $Q^d_p$ \cite{arXiv:2510.13348}, it seems natural to ask the following:
\begin{question}
    Under which vertex and edge expansion conditions on the host graph $G$ and under which condition on $d$ and $p$, is it true that \textbf{whp} the giant component $L_1$ in $G_p$ satisfies
    \[
        \mathrm{diam}(L_1) = O(d + \log n)?
    \]
\end{question}
In Theorem~\ref{general_large_diam} we also observe a quite large conceptual gap between our families of graphs with good vertex and edge expansion from one side, and well-studied graphs such as random subgraphs of $K_n$ and $(n,d,\lambda)$-graphs, from the other side. It is therefore natural to ask what makes the difference:
\begin{question}
    What further assumptions on the host graph $G$, beyond those in Theorems \ref{isop},\ref{ERCP1} and \ref{L1}, are needed to obtain that \textbf{whp} the giant component $L_1$ in $G_p$ satisfies
    \[
        \mathrm{diam}(L_1) = \Theta(\log n)?
    \]
\end{question}
Observe that in order to achieve the upper bound $\mathrm{diam}(L_1) = O(\log n)$, it suffices that the giant component $L_1$ in $G_p$ satisfies the following \textbf{whp}: $|N_{L_1}(U)|=\Omega(|U|)$ for any connected $U\subseteq L_1$ with $|U|=\Omega(\log n)$ (see, for example, the proof of Corollary \ref{positive_diam}). The lower bound is satisfied by Theorem~\ref{like_Chung_Lu}.

Another direction for future research is to investigate the corresponding results for site percolation. In \cite{arXiv:2603.19038}, it was shown that the typical emergence of a giant component under site percolation is guaranteed under certain expansion assumptions on the host graph. However, it was also demonstrated that, in this respect, site and bond percolation exhibit somewhat different behavior. It would therefore be interesting to investigate whether further differences arise in the combinatorial properties of the giant component under site percolation.

We finish this section with a conjecture. We have shown that under standard expansion properties of the host graph we have a typical vertex expansion for linear-sized sets and for connected sets in the giant component. This typical vertex expansion is by a factor of $\frac{1}{d\log d}$. We believe that the following should also be true:

\begin{conjecture}\label{conj}
    Let $G$ be a $d$-regular $n$-vertex graph that satisfies the assumptions of Theorem~\ref{ERCP1} and let $p=\frac{1+\epsilon}{d}$. Then, for every constant $\delta>0$, there exists a constant $\beta\coloneqq \beta(\epsilon,c_1,c_2,c_3,C,\alpha,\delta)>0$ such that \textbf{whp} every connected set $S\subseteq V(L_1)$ of the giant component $L_1\subseteq G_p$ of size $\delta n\le |S|\le \frac{|V(L_1)|}{2}$ satisfies
    \[
        |N_{G_p}(S)|\ge \frac{\beta}{d}|S|.
    \]
\end{conjecture}

Two follow-up questions regarding this conjecture may be natural here: first, does the lower bound $\delta n$ have to be linear, or can Conjecture~\ref{conj} be established for sublinear sets as well? Furthermore, it is interesting to ask whether the connectivity assumption in Conjecture~\ref{conj} can be removed, perhaps by choosing an appropriate $\delta$ that depends on the size of the giant component.

One of the motivations for Conjecture~\ref{conj} is a related result for infinite connected transitive graphs. Essentially, Diskin, Lichev, Ramanan Radhakrishnan, Sudakov, and Verhoeven \cite{arXiv:2610.38941} showed that if an infinite graph has an isoperimetric profile $k\rightarrow f(k)$ (see the definition in \cite{arXiv:2610.38941}), then almost surely the infinite cluster in the percolated graph has an anchored isoperimetric profile $k\rightarrow c\cdot f(k)$, where $c=c(p)>0$ is a constant.

\section*{Acknowledgments} 
Part of this work was carried out during visits by the second and third authors to the Erwin Schrödinger International Institute for Mathematics and Physics (ESI) in Vienna in October 2026. They gratefully acknowledge ESI’s support. 
This research was supported in part by the Austrian Science Fund (FWF) [10.55776/F1002] (M. Kang) and NSF-BSF grant 2023688 (M. Krivelevich). For open access purposes, the authors have applied a CC BY public copyright license to any author accepted manuscript version arising from this submission.

\paragraph{AI Disclosure:} AI tools were used to assist with proofreading and searching of existing literature.

\bibliographystyle{abbrv}
\bibliography{Bib}
\newpage

\appendix
\section{Properties of $G$ from Theorem~\ref{no_L1}}\label{app: 5}
In Section \ref{pno_L1} we aimed to construct a graph $G$ that satisfies $\mathrm{iso}(b)$ and $\mathrm{EE}(Cd,\alpha)$, and, moreover, \textbf{whp} all components of $G_p$ are of size at most $\frac{9}{\alpha}Cd\log n$. In this section we show that the graph constructed in Section \ref{pno_L1} indeed satisfies properties $\mathrm{iso}(b)$ and $\mathrm{EE}(Cd,\alpha)$.

\begin{claim}\label{pprop2S}
    $G$ satisfies $\mathrm{EE}(Cd,\alpha)$.
\end{claim}

\begin{proof}
    Let $U\subseteq V(G)$ such that $|U|\le Cd$. Define $U_0=U\cap V(H_0)$ and for every $1\le i\le n_1$, we set $U_i=U\cap V(H_i)$ and $b_i\coloneqq |U_i\cap F_i|$. Observe that each $U_i$ touches at most $b_i$ edges from $E_i$.
    Notice that 
    \begin{align*}
        e_G(U,U^c)=e_{H_0}(U_0,U_0^c)&+\sum_{i: v_i\notin U_0}\left(e_{H_i}(U_i,U_i^c)+e_G(\{v_i\}, U_i\cap F_i)\right)\\
        &+\sum_{i: v_i\in U_0}\left(e_{H_i}(U_i,U_i^c)+e_G(\{v_i\}, F_i\backslash U_i)\right).
    \end{align*}
    We analyze each term separately.
    By Claim~\ref{eml},
    \begin{align*}
        e_{H_0}(U_0,U_0^c)\ge \frac{d_1-\lambda}{n_1}|U_0|(n_1-|U_0|)\ge \left(1-\frac{\lambda}{d_1}\right)\left(1-\frac{Cd}{n_1}\right)d_1|U_0|\ge \left(1-\frac{\alpha}{2}\right)d_1|U_0|,
    \end{align*}
    where the last step holds since $\lambda=O\left(\sqrt{d_1}\right)$, $d_1=\omega(1)$ and $n=\omega\left(d^2\right)$.
    Let $1\le i\le n_1$, and assume that $v_i\notin U_0$. Then, $e_G(\{v_i\}, U_i\cap F_i)=b_i$ and since $H_i$ is a copy of $\left(\frac{4}{\alpha}Cd,d,\mu\right)$-graph minus a matching of size $\frac{d-d_1}{2}$, then by Claim~\ref{eml}:
    \begin{align*}
        e_{H_i}(U_i,U_i^c)\ge \frac{d-\mu}{\frac{4}{\alpha}Cd}|U_i|\left(\frac{4}{\alpha}Cd-|U_i|\right)-b_i.
    \end{align*}
    Together we yield
    \begin{align*}
        e_G(\{v_i\}, U_i\cap F_i)+e_{H_i}(U_i,U_i^c)&\ge \frac{d-\mu}{\frac{4}{\alpha}Cd}|U_i|\left(\frac{4}{\alpha}Cd-|U_i|\right)\\
        &\ge \left(1-\frac{\mu}{d}\right)\left(1-\frac{\alpha}{4}\right)d|U_i|\ge \left(1-\frac{\alpha}{2}\right)d|U_i|,
    \end{align*}
    where the last step holds since $\mu=O\left(\sqrt{d}\right)$, $d=\omega(1)$.

    Now turn to the case where $v_i\in U_0$. Similarly to the previous case, we have
    \begin{align*}
        e_{H_i}(U_i,U_i^c)+e_G(\{v_i\}, F_i\backslash U_i)&\ge \frac{d-\mu}{\frac{4}{\alpha}Cd}|U_i|\left(\frac{4}{\alpha}Cd-|U_i|\right)-b_i+d-d_1-b_i\\
        &\ge \left(1-\frac{\alpha}{2}\right)d|U_i| +(d-d_1)-2b_i.
    \end{align*}
    Notice that 
    \begin{align*}
        \sum_{i: v_i\in U_0}\left(e_{H_i}(U_i,U_i^c)+e_G(\{v_i\}, F_i\backslash U_i)\right)&\ge \sum_{i: v_i\in U_0}\left[\left(1-\frac{\alpha}{2}\right)d|U_i| +(d-d_1)-2b_i\right]\\
        &\ge \sum_{i: v_i\in U_0} \left(1-\frac{\alpha}{2}\right)d|U_i|+(d-d_1)|U_0|-2\sum_{i: v_i\in U_0}b_i\\
        &\ge \sum_{i: v_i\in U_0}\left(1-\frac{\alpha}{2}\right)d|U_i|+\left(1-\frac{\alpha}{2}\right)(d-d_1)|U_0|-2|U|.
    \end{align*}
    Altogether, we get
    \begin{align*}
        e_G(U,U^c)&\ge \left(1-\frac{\alpha}{2}\right)d_1|U_0|+\left(1-\frac{\alpha}{2}\right)(d-d_1)|U_0|-2|U|+\sum_{i=1}^{n_1} \left(1-\frac{\alpha}{2}\right)d|U_i|\\
        &= \left(1-\frac{\alpha}{2}\right)d|U_0|-2|U|+\sum_{i=1}^{n_1} \left(1-\frac{\alpha}{2}\right)d|U_i|\\
        &= \left(1-\frac{\alpha}{2}\right)d|U|-2|U|\ge (1-\alpha)d|U|,
    \end{align*}
    where the last inequality holds since $d=\omega(1)$.
\end{proof}

\begin{claim}
    $G$ satisfies $\mathrm{iso}(b)$, where $b\coloneqq \frac{\alpha}{60\sigma e^2 C}$.
\end{claim}

\begin{proof}
    Let $U\subseteq V(G)$, with $|U|\le \frac{n}{2}$. We show that $\frac{e_G(U,U^c)}{|U|}\ge b$. Set $U_0=U\cap V(H_0)$ and for every $1\le i\le n_1$, $U_i=U\cap (\{v_i\}\cup V(H_i))$. 
    Suppose first that $|U_0|\ge \frac{4|U|}{5(\frac{4}{\alpha}Cd+1)}$ and $|U_0|\le \frac{|V(H_0)|}{2}$. Notice that
    \begin{align*}
        \frac{e_G(U,U^c)}{|U|}\ge \frac{e_{H_0}(U_0,V(H_0)\backslash U_0)}{|U|}.
    \end{align*}
    Now, since $H_0$ is a $(n_1,d_1,\lambda)$-graph with $\lambda=O\left(\sqrt{d_1}\right)$, then by Claim~\ref{eml}, for every $W\subseteq V(H_0)$ with $|W|\le \frac{|V(H_0)|}{2}$ we have $e_{H_0}(W,W^c)\ge \frac{d_1}{3}|W|$, and hence:
    \begin{align*}
        \frac{e_{H_0}(U_0,V(H_0)\backslash U_0)}{|U|}\ge \frac{\frac{d_1|U_0|}{3}}{|U|}\ge \frac{4d_1}{15\left(\frac{4}{\alpha}Cd+1\right)}\ge \frac{4}{15\sigma e^2}\cdot \frac{d}{\frac{4}{\alpha}Cd+1}\ge \frac{\alpha}{60\sigma e^2 C}=b,
    \end{align*}
    where we used the facts that $d_1=\frac{d}{\sigma e^2}$ and $b=\frac{\alpha}{60\sigma e^2 C}$.
    
    Now, suppose that $|U_0|\ge \frac{4|U|}{5(\frac{4}{\alpha}Cd+1)}$ and $\frac{|V(H_0)|}{2}\le |U_0|\le \frac{2|V(H_0)|}{3}$. Thus, we have
    \begin{align*}
        \frac{e_G(U,U^c)}{|U|}&\ge \frac{e_{H_0}(U_0,V(H_0)\backslash U_0)}{|U|}\ge \frac{d_1|V(H_0)\backslash U_0|}{3|U|}\ge \frac{d_1|V(H_0)|}{9|U|}\\
        &\ge \frac{d_1|V(H_0)|}{9\cdot\frac{5\left(\frac{4}{\alpha}Cd+1\right)|V(H_0)|}{6}}=\frac{2}{15}\cdot \frac{d_1}{\frac{4}{\alpha}Cd+1}\ge b.
    \end{align*}
    
    We turn to the case where $|U_0|< \frac{4|U|}{5(\frac{4}{\alpha}Cd+1)}$. Observe that
    \begin{align}\label{aux2}
        \frac{4|U|}{5}>|U_0|\left(\frac{4}{\alpha}Cd+1\right)= \left(\sum_{i:v_i\in U_i}1\right)\left(\frac{4}{\alpha}Cd+1\right)= \bigg|\bigcup_{i:v_i\in U_i}\left(V(H_i)\cup \{v_i\}\right)\bigg|.
    \end{align} 
    Since we can write
    \begin{align*}
        U=\left(\bigcup_{i:v_i\in U_i} U_i\right) \cup \left(\bigcup_{i:v_i\notin U_i\newline } U_i\right)\subseteq \left(\bigcup_{i:v_i\in U_i} \left(V(H_i)\cup \{v_i\}\right)\right) \cup \left(\bigcup_{i:v_i\notin U_i} U_i\right),
    \end{align*}
    then by \eqref{aux2}, we obtain
    \begin{align*}
        \bigg|\bigcup_{i:v_i\notin U_i} U_i\bigg|\ge |U|-\bigg|\bigcup_{i:v_i\in U_i} \left(V(H_i)\cup \{v_i\}\right)\bigg|\ge \frac{|U|}{5}.
    \end{align*}
    Now, as $U=\bigcup_{i=1}^{n_1} U_i$, we have
    \begin{align}\label{pain2}
        \frac{e_G(U,U^c)}{|U|}\ge \frac{1}{5}\min_{i:v_i\notin U_i} \frac{e_G(U_i,\{v_i\}\cup V(H_i)\backslash U_i)}{|U_i|}.
    \end{align}
    Set $1\le i\le n_1$ such that $v_i\notin U_i$. Notice that if $U_i=V(H_i)$, then 
    \begin{align*}
        \frac{e_G(U_i,\{v_i\}\cup V(H_i)\backslash U_i)}{5|U_i|}=\frac{e_G(U_i,\{v_i\})}{5|U_i|}=\frac{d-d_1}{5\cdot \frac{4}{\alpha}Cd}\ge \frac{\alpha}{60C}\ge b,
    \end{align*}
    hence, we can assume that $|U_i|,|V(H_i)\backslash U_i|\ge 1$.
    Similarly to the proof of Claim~\ref{pprop2S}, set $b_i=\left|U_i\cap F_i\right|$. Observe that
    \begin{align*}
        e_{H_i}(U_i,V(H_i)\backslash U_i)\ge \frac{d-\mu}{\frac{4}{\alpha}Cd} |U_i||V(H_i)\backslash U_i|-|U_i\cap F_i|= \frac{d-\mu}{\frac{4}{\alpha}Cd} |U_i||V(H_i)\backslash U_i|-b_i.
    \end{align*}
    Therefore, we have
    \begin{align*}
        e_G(U_i,\{v_i\}\cup V(H_i)\backslash U_i)&= e_{H_i}(U_i,V(H_i)\backslash U_i)+e_G(U_i,\{v_i\})\ge \frac{d-\mu}{\frac{4}{\alpha}Cd} |U_i||V(H_i)\backslash U_i|-b_i+b_i\\
        &=\frac{d-\mu}{\frac{4}{\alpha}Cd} |U_i||V(H_i)\backslash U_i|\ge \frac{d-\mu}{2}\ge \frac{d}{4},
    \end{align*}
    where the penultimate inequality holds since $|U_i|,|V(H_i)\backslash U_i|\ge 1$ and since $|U_i|+|V(H_i)\backslash U_i|=\frac{4}{\alpha}Cd$. The last inequality holds since $\mu=O\left(\sqrt{d}\right)$.
    Now,
    \begin{align*}
        \frac{e_G(U_i,\{v_i\}\cup V(H_i)\backslash U_i)}{5|U_i|}\ge \frac{d}{20|U_i|}\ge \frac{d}{20\cdot\frac{4}{\alpha}Cd}\ge \frac{\alpha}{60\sigma e^2 C}=b.
    \end{align*}
    Overall, by \eqref{pain2}, we have  
    \begin{align*}
        \frac{e_G(U,U^c)}{|U|}\ge b.
    \end{align*}
    We have only one more case to consider, which is when $|U_0|\ge \frac{2|V(H_0)|}{3}$. In this case, there exists $U'\subseteq U_0$ with $|U'|\ge \frac{|V(H_0)|}{6}$ such that for every $v_i\in U'$, $U_i\backslash \{v_i\}\subsetneq V(H_i)$. Let $1\le i\le n_1$, such that $v_i\in U'$. If $U_i\backslash \{v_i\}=\emptyset$, then 
    \begin{align*}
        e_G(U_i,V(H_i)\backslash U_i)=e_G(\{v_i\},V(H_i))=d-d_1\ge \frac{d}{4}.
    \end{align*}
    We now assume that $U_i\backslash \{v_i\}\neq \emptyset$. Denote by $U_i'\coloneqq U_i\backslash \{v_i\}$. Hence, $|U_i'|,|V(H_i)\backslash U_i'|\ge 1$. If $1\le|U_i'|\le 2$ or $1\le|V(H_i)\backslash U_i'|\le 2$, then
    \begin{align*}
        e_G(U_i,V(H_i)\backslash U_i)\ge e_{H_i}(U_i',V(H_i)\backslash U_i')\ge \frac{d}{4}.
    \end{align*}
    Now assume that $|U_i'|,|V(H_i)\backslash U_i'|\ge 3$. Then
    \begin{align*}
        e_G(U_i,V(H_i)\backslash U_i)&\ge e_{H_i}(U_i',V(H_i)\backslash U_i')\ge \frac{d-\mu}{\frac{4}{\alpha}Cd} |U_i'||V(H_i)\backslash U_i'|-b_i\ge \frac{d-\mu}{2}\cdot 3-b_i\\
        &\ge \frac{3(d-\mu)}{2}-(d-d_1)\ge  \frac{d}{4}.
    \end{align*}
    Overall, we have
    \begin{align*}
        e_G(U,U^c)\ge \sum_{i:v_i\in U'}e_G(U_i,V(H_i)\backslash U_i)\ge \sum_{i:v_i\in U'} \frac{d}{4}\ge \frac{d|V(H_0)|}{24}.
    \end{align*}
    Therefore,
    $\frac{e_G(U,U^c)}{|U|}\ge \frac{d|V(H_0)|}{24\cdot \frac{n}{2}}=\frac{d\cdot \frac{n}{\frac{4}{\alpha}Cd+1}}{12n}\ge b$,
    as desired.
\end{proof}

\section{Properties of $G$ from Theorem~\ref{no_expansion}}\label{app: 6}
In this section we prove that the graph $G$, constructed in Section \ref{Pno_expansion}, satisfies properties $\mathrm{iso}(c_1)$, $\mathrm{VE}(c_2\log n,\alpha)$ and $\mathrm{EE}(c_2 d\log n,\alpha)$.

\begin{claim}
    $G$ satisfies $\mathrm{iso}(c_1)$.
\end{claim}

\begin{proof}
    Let $U\subseteq V(G)$ be with $|U|\le \frac{n}{2}$. Set $U_0=V(H^0)\cap U$. We prove the claim in the case when $|U_0|\le \frac{n}{4}$, as the proof of the complementary case is very similar. As long as $|U_0|\ge \frac{|U|}{6}$, using Claim~\ref{eml} we have 
    \begin{align*}
        e_G(U,U^c)\ge e_{H^0}(U_0,U_0^c)\ge \frac{d-C-\tau}{\frac{n}{2}}\cdot |U_0||U_0^c|\ge \frac{d-C-\tau}{2}\cdot |U_0|\ge \frac{d}{4}\cdot |U_0|\ge \frac{d}{24}|U|.
    \end{align*}
    Now assume that $|U_0|\le \frac{|U|}{6}$. For every $1\le i\le k$, set $S_i=A_i\cap U$ and $T_i=B_i\cap U$ and also denote $x_i\coloneqq |S_i|$ and $y_i\coloneqq |T_i|$. Denote by $I\subseteq [k]$ the set of indices for which $x_i+y_i\le 1.5d\log n$. We claim that 
    \renewcommand\qedsymbol{$\blacksquare$}
    \begin{claim}\label{aux_claim}
        \begin{align*}
            \sum_{i\in I}(x_i+y_i)>\sum_{i\notin I}(x_i+y_i).
        \end{align*}
    \end{claim}
    \begin{proof}
        Assume otherwise. Then we have
        \begin{align*}
            \sum_{i\in I}(x_i+y_i)\le \sum_{i\notin I}(x_i+y_i)\le 3\sum_{i\notin I}y_i,
        \end{align*}
        where the last inequality holds by the following observation: notice that for $i\notin I$ we have $x_i+y_i> 1.5d\log n$, and since $x_i\le |A_i|=d\log n$, we yield that $x_i\le 2y_i$. Thus, 
        \begin{align*}
            \sum_{i=1}^k (x_i+y_i)\le 6\sum_{i\notin I}y_i,
        \end{align*} 
        and therefore:
        \begin{align*}
            \frac{|U|}{6}=\frac{1}{6}\sum_{i=1}^k (x_i+y_i)\le \sum_{i\notin I}y_i\le |U_0|<\frac{|U|}{6},
        \end{align*}
        which is a contradiction.        
    \end{proof}
    \renewcommand\qedsymbol{$\square$}
    \noindent Now, by Claim~\ref{eml} for the graphs $\{G[A_i,B_i]\}_{i\in[k]}$, for every $i\in I$, we have that
    \[
        e_{G[A_i,B_i]}(S_i\cup T_i,(S_i\cup T_i)^c)\ge \frac{C-\lambda}{2d\log n}(x_i+y_i)\left(2d\log n-x_i-y_i\right).
    \]
    Hence, by Claim~\ref{aux_claim}, we conclude that
    \begin{align*}
        e_G(U,U^c)&\ge \sum_{i\in I}e_{G[A_i,B_i]}(S_i\cup T_i,(S_i\cup T_i)^c)\ge \sum_{i\in I} \frac{C-\lambda}{2d\log n}(x_i+y_i)\left(2d\log n-x_i-y_i\right)\\
        &\ge \sum_{i\in I}\frac{C-\lambda}{4}(x_i+y_i)\ge \frac{C}{8}\sum_{i\in I}(x_i+y_i)\ge \frac{C}{16}\sum_{i=1}^k(x_i+y_i)\\
        &=\frac{C}{16}|U|=\frac{1}{160\epsilon}|U|\ge c_1|U|,
    \end{align*}
    where the last inequality holds since $c_1<\frac{1}{160\epsilon}$.
\end{proof}

Next, to establish vertex expansion, we need the following result, which is a weaker version of Lemma 16 in \cite{arXiv:2602.17567} (but suffices for our purposes).

\begin{lemma}\label{G_n,d vertex expansion}
    Let $c>0$ be a sufficiently small constant. Let $H$ be a random $d$-regular graph on $n$ vertices, such that $d=\omega_n(1)$ and $d\le \frac{n}{2}$. Then, the following holds \textbf{whp}: every subset $U\subseteq V(H)$ of size $|U|\le c\frac{n}{d}$ satisfies
    \[
        |N_H(U)|\ge (1-200c)d|U|.
    \]
\end{lemma}

\begin{claim}
    The graph $G$ satisfies $\mathrm{VE}(c_2\log n,\alpha)$.
\end{claim}

\begin{proof}
    Let $U\subseteq V(G)$ be a set of vertices of size $|U|\le c_2\cdot \log n$. For every $0\le i\le k$, define $U_i=V(H^i)\cap U$.
    As $|U|\le c_2\cdot \log n$, we have $|U_i|\le c_2\cdot \log n$, for every $0\le i\le k$.
    
    The graph $H^i$ is a copy of the random $(d-C)$-regular graph on $d\log n$ vertices ($1\le i\le k$) and $H^0$ is a copy of the random $(d-C)$-regular graph on $\frac{n}{2}$ vertices. Applying Lemma~\ref{G_n,d vertex expansion} on $H^i$ for every $0\le i\le k$ (with $c\coloneqq \frac{\alpha}{400}$ and $c_2=c$), we obtain that
    \begin{align*}
        |N_{H^i}(U_i)|\ge \left(1-\frac{\alpha}{2}\right)(d-C)|U_i|\ge (1-\alpha)d|U_i|,
    \end{align*}
    where the last inequality holds since $d=\omega(1)$. Therefore, 
    \[
        |N_G(U)|\ge \sum_{i=0}^{k} |N_{H^i}(U_i)|\ge (1-\alpha)d\cdot \sum_{i=0}^{k} |U_i| \ge (1-\alpha)d|U|.\qedhere
    \] 
\end{proof}

Finally, we claim the desired edge expansion.

\begin{claim}
    The graph $G$ satisfies $\mathrm{EE}\left(\frac{\alpha}{2}\cdot d\log n,\alpha\right)$.
\end{claim}

\begin{proof}
    Let $U\subseteq V(G)$ be a set of vertices of size $|U|\le \frac{\alpha}{2}\cdot d\log n$. For every $0\le i\le k$, define $U_i=V(H^i)\cap U$. Since $|U|\le \frac{\alpha}{2}\cdot d\log n$, we have $|U_i|\le \frac{\alpha}{2}\cdot d\log n$, for every $0\le i\le k$. Hence, by Claim~\ref{eml}, for every $1\le i\le k$, we have
    \begin{align*}
        e_{H^i}(U_i,U_i^c)\ge \frac{(d-C)-\nu}{d\log n}|U_i||U_i^c|\ge \left(1-\frac{\alpha}{2}\right)(d-C-\nu)|U_i|\ge (1-\alpha)d|U_i|,
    \end{align*}
    where the last inequality holds since $\nu=O\left(\sqrt{d-C}\right)$ and $d=\omega(1)$. Similarly, $e_{H^0}(U_0,U_0^c)\ge (1-\alpha)d|U_0|$. Therefore,
    \[
        e_G(U,U^c)\ge \sum_{i=0}^{k} e_{H^i}(U_i,U_i^c) \ge \sum_{i=0}^{k}(1-\alpha)d|U_i|=(1-\alpha)d|U|.\qedhere
    \]
\end{proof}

\section{Properties of $G$ from Theorem~\ref{general_large_diam}}\label{app: 7}
In this section we prove that the graph $G$, constructed in Section \ref{P_larege_diam}, satisfies properties $\mathrm{iso}(c_1)$, $\mathrm{VE}(c_2\cdot \frac{t}{d},\alpha)$ and $\mathrm{EE}(\alpha t,\alpha)$.

\begin{claim}\label{pprop0U}
    $G$ satisfies $\mathrm{iso}(c_1)$.
\end{claim}

\begin{proof}
    Consider $U\subseteq V(G)$ such that $|U|\le \frac{n}{2}$. For every $1\le i\le \frac{n}{t}$, set $U_i\coloneqq U\cap V(H_i)$ and $a_i\coloneqq |U_i|$. Observe that by applying Claim~\ref{eml} on $H_i$ (for every $1\le i\le \frac{n}{t}$), we have
    \begin{align*}
        e_G(U,U^c)&=\sum_{i=1}^{n/t} e_{H_{i}}(U_i,U_i^c)\ge \sum_{i=1}^{n/t} \frac{\delta d-\lambda_2}{2t}\cdot a_i(2t-a_i)\ge \frac{\delta d}{3t} \sum_{i=1}^{n/t}\left(2ta_i-a_i^2\right)\\
        &=\frac{2\delta d}{3}\sum_{i=1}^{n/t}a_i-\frac{\delta d}{3t}\sum_{i=1}^{n/t}a_i^2=\frac{4\delta d}{3}|U|-\frac{\delta d}{3t}\sum_{i=1}^{n/t}a_i^2,
    \end{align*}
    where the second inequality holds since $\lambda_2= O\left(\sqrt{d}\right)$ and since $d=\omega(1)$, and the last equality holds since $\sum_{i=1}^{n/t} a_i=2|U|$. 
    
    Now, let us first assume that $|U|\ge t$.
    Using the fact that $\sum_{i=1}^{n/t} a_i=2|U|$, the maximal value of $\sum_{i=1}^{n/t} a_i^2$ is obtained when one of them is equal to $t$, $\lfloor\frac{|U|}{t}\rfloor-1$ of $\{a_i\}_{i=1}^{n/t}$ are equal to $2t$, one of them is equal to $t+t\cdot \left\{\frac{|U|}{t}\right\}$, one of them is equal to $t\cdot \left\{\frac{|U|}{t}\right\}$ and the rest are equal to $0$. Hence,
    \begin{align*}
        \sum_{i=1}^{n/t}a_i^2&\le t^2+\left(\frac{|U|}{t}-\left\{\frac{|U|}{t}\right\}-1\right)\cdot (2t)^2+\left(t+t\cdot\left\{\frac{|U|}{t}\right\}\right)^2+\left(t\cdot\left\{\frac{|U|}{t}\right\}\right)^2\\
        &=4t|U|-2t^2-2t^2\left\{\frac{|U|}{t}\right\}\left(1-\left\{\frac{|U|}{t}\right\}\right)\le 4t|U|-2t^2.
    \end{align*}
    In conclusion, we have
    \begin{align*}
        e_G(U,U^c)&\ge \frac{4\delta d}{3}|U|-\frac{\delta d}{3t}\sum_{i=1}^{n/t}a_i^2\ge \frac{4\delta d}{3}|U|-\frac{\delta d}{3t}\cdot \left(4t|U|-2t^2\right)=\frac{2}{3}\delta dt\ge \frac{2}{3}c_1n\ge c_1\frac{n}{2}\ge c_1|U|,
    \end{align*}
    where we use here the fact that $t\ge \frac{c_1}{\delta}\cdot \frac{n}{d}$.

    Now assume that $|U|<t$. The maximal value of $\sum_{i=1}^{n/t} a_i^2$ is obtained when exactly two of the $\{a_i\}_{i=1}^{n/t}$ are equal to $|U|$ and the rest are equal to 0. Therefore,
    \[
        e_G(U,U^c)\ge \frac{4\delta d}{3}|U|-\frac{\delta d}{3t}\sum_{i=1}^{n/t}a_i^2\ge \frac{4\delta d}{3}|U|-\frac{2\delta d}{3t}|U|^2\ge \frac{2d}{3}|U|\ge c_1|U|,
    \]
    where the penultimate inequality holds since $|U|<t$.
\end{proof}

\begin{claim}
    $G$ satisfies $\mathrm{VE}(c_2\cdot \frac{t}{d},\alpha)$.
\end{claim}

\begin{proof}
    Let $U\subseteq V(G)$ with $|U|\le c_2\frac{t}{d}$. For every $1\le i\le \frac{n}{t}$, set $T_i=U\cap S_i$. We have $|T_i|\le c_2\frac{t}{d}$ for every $1\le i\le \frac{n}{t}$. Observe that
    \[
        |N_G(U)|\ge\sum_{i=1}^{n/t} |(N_G(T_i)\cap S_i)|.
    \]
    The graph $F_i=G[S_i]$ is a copy of the random $\left(1-\frac{\alpha}{2}\right)d$-regular graph on $t$ vertices for every $1\le i\le \frac{n}{t}$. Hence, by Lemma~\ref{G_n,d vertex expansion} (with $c\coloneqq \frac{\alpha}{400}$ and $c_2=c$), we obtain that for every $0\le i\le \frac{n}{t}$ 
    \[
        |(N_G(T_i)\cap S_i)|\ge \left(1-\frac{\alpha}{2}\right)\left(1-\frac{\alpha}{2}\right)d|T_i|\ge (1-\alpha)d|T_i|.
    \]
    Hence,
    \[
        |N_G(U)|\ge\sum_{i=1}^{n/t} |(N_G(T_i)\cap S_i)|\ge (1-\alpha)d\sum_{i=1}^{n/t}|T_i|=(1-\alpha)d|U|.\qedhere
    \]
\end{proof}

\begin{claim}
    $G$ satisfies $\mathrm{EE}(\frac{\alpha}{2} t,\alpha)$.
\end{claim}

\begin{proof}
    Let $U\subseteq V(G)$ with $|U|\le \frac{\alpha}{2} t$. For every $1\le i\le \frac{n}{t}$, set $T_i\coloneqq U\cap S_i$ and $a_i\coloneqq |T_i|$.
    Notice that
    \[
        e_G(U,U^c)\ge\sum_{i=1}^{n/t} e_{F_i}(T_i,T_i^c).
    \]
    Recall that for every $1\le i\le \frac{n}{t}$, the graph $F_i=G[S_i]$ is a $\left(t,\left(1-\frac{\alpha}{2}\right)d,\tau\right)$-graph with $\tau=O\left(\sqrt{d}\right)$. Applying Claim~\ref{eml} on $F_i$ (for every $1\le i\le \frac{n}{t}$), we obtain
    \[
        e_{F_i}(T_i,T_i^c)\ge \frac{\left(1-\frac{\alpha}{2}\right)d-\tau}{t}\cdot a_i(t-a_i)\ge \left(1-\frac{\alpha}{2}-\frac{\tau}{d}\right)\left(1-\frac{\alpha}{2}\right)d a_i\ge (1-\alpha)da_i,
    \]
    where in the penultimate inequality we used $a_i\le \frac{\alpha}{2} t$ and in the last inequality, we used the facts that $\tau= O\left(\sqrt{d}\right)$ and that $d=\omega(1)$. To conclude, 
    \[
        e_G(U,U^c)\ge\sum_{i=1}^{n/t} e_{F_i}(T_i,T_i^c)\ge (1-\alpha)d\sum_{i=1}^{n/t}a_i=(1-\alpha)d|U|.
    \]
\end{proof}
\end{document}